\documentclass[10pt]{amsart}

\usepackage[utf8]{inputenc} 
\usepackage[T1]{fontenc}
\usepackage{url}
\usepackage[margin=1.3in]{geometry}
\usepackage{parskip}

\usepackage{amssymb} 
\usepackage{amsmath}
\usepackage{amsthm} 
\usepackage{mathabx} 
\usepackage{mathtools}

\allowdisplaybreaks 

\usepackage{enumerate}
\usepackage[shortlabels]{enumitem}

\newcommand{\Q}{\mathbb{Q}}
\newcommand{\Z}{\mathbb{Z}}
\newcommand{\del}{\partial}

\usepackage{graphicx} 
\usepackage{tikz}
\usepackage[center]{caption} 

\usepackage[hidelinks,bookmarksopen, bookmarksopenlevel=0]{hyperref}
\usepackage{bookmark}

\newtheorem{theor}{Theorem} 
\newtheorem{thm}{Theorem}[section]
\newtheorem{prop}[thm]{Proposition}

\newtheorem{lemma}[thm]{Lemma}
\newtheorem{cor}[thm]{Corollary}
\newtheorem*{thm*}{Theorem}
\newtheorem*{prop*}{Proposition}
\newtheorem*{cor*}{Corollary}
\newtheorem*{prop:sxpiece}{Proposition~\ref{prop:sxpiece}}
\newtheorem*{thm:composite}{Theorem~\ref{thm:composite}}
\newtheorem*{thm:cableSFS}{Theorem~\ref{thm:cableSFS}}

\theoremstyle{definition}
\newtheorem{defn}[thm]{Definition}

\newtheorem{rem}[thm]{Remark}

\newtheorem*{acknow}{Acknowledgements}

\theoremstyle{remark}

\begin{document}

\title{Characterizing slopes for satellite knots}
\author{Patricia Sorya}
\address{Département de mathématiques, Université du Québec à Montréal}
\email{sorya.patricia@courrier.uqam.ca}

\thanks{Supported by FRQNT doctoral reasearch grant 305903}

\begin{abstract}
A slope $p/q$ is said to be characterizing for a knot $K$ if the homeomorphism type of the $p/q$-Dehn surgery along $K$ determines the knot up to isotopy. Extending previous work of Lackenby and McCoy on hyperbolic and torus knots respectively, we study satellite knots to show that for a knot $K$, any slope $p/q$ is characterizing provided $|q|$ is sufficiently large. In particular, we establish that every non-integral slope is characterizing for a composite knot. Our approach consists of a detailed examination of the JSJ decomposition of a surgery along a knot, combined with results from other authors giving constraints on surgery slopes that yield manifolds containing certain surfaces. 
\end{abstract}
\maketitle


  \section{Introduction}

A non-trivial slope $p/q$ is said to be \emph{characterizing} for a knot $K$ in $S^3$ if whenever there exists an orientation-preserving homeomorphism $S^3_{K}(p/q) \cong S^3_{K'}(p/q)$ between the $p/q$-Dehn surgery along $K$ and the $p/q$-Dehn surgery along some knot $K'$, then $K = K'$, where ``$=$'' denotes an equivalence of knots up to isotopy.
In \cite{lack}, Lackenby proved that every knot has infinitely many characterizing slopes by showing that any slope is characterizing for a knot $K$, provided $|p|\leq|q|$ and $|q|$ is sufficiently large.

The main theorem of the present paper strengthens this result.
\begin{theor}\label{thm:main}
Let $K$ be a knot in $S^3$. Then any slope $p/q$ is characterizing for $K$, provided $|q|$ is sufficiently large.
\end{theor}

In \cite{kronheimer_monopoles_2007},  Kronheimer, Mrowka, Ozsváth and Szabó proved that all non-trivial slopes are characterizing for the unknot. McCoy showed in \cite{mcc} that if $K$ is a torus knot, there are only finitely many non-integral slopes that are non-characterizing for $K$, thus giving the torus knot case of the theorem. Lackenby showed the hyperbolic knot case in \cite{lack}. In this paper, we establish the theorem for any knot by studying the case of satellite knots. 

The extension to satellite knots requires a distinct approach, as it cannot be simply derived from the cases of hyperbolic and torus knots. This is due to the presence of essential tori in the exterior of a satellite knot, which lead to a non-trivial JSJ decomposition of the knot's exterior. Hence, Dehn surgery along a satellite knot involves attaching a solid torus to a torus boundary component of a manifold that is not a knot exterior. Our strategy therefore consists of an in-depth analysis of the topology of Dehn fillings of manifolds that arise as JSJ pieces of a knot exterior, along with a description of the gluing between these manifolds through the distance between specific slopes. In particular, we rely on the rigidity of Seifert fibred structures, as well as results pertaining to fillings of non-Seifert fibred manifolds that contain certain surfaces.

Moreover, the ideas employed in the proof of Theorem \ref{thm:main} can be adapted to derive explicit bounds on $|q|$ for some families of satellite knots. We obtain the following result for composite knots.

\begin{theor}\label{thm:composite}
If $K$ is a composite knot, then every non-integral slope is characterizing for $K$.
\end{theor}

\begin{figure}[h]
    \centering\begin{tikzpicture}
    \node at (0,0) {\includegraphics[scale=0.33]{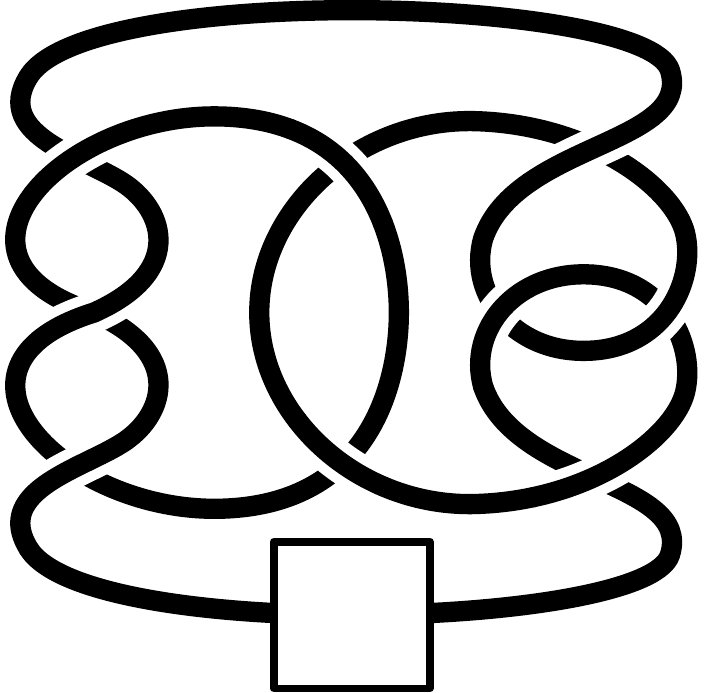}};
    \node at (0,-43pt) {$K$};
    \end{tikzpicture}
    \caption{The set of characterizing slopes for the connected sum of $9_{42}$ and any non-trivial knot $K$ is $\Q \setminus \Z$}
    \label{fig:bakermotegi}
\end{figure}

Baker and Motegi constructed composite knots for which every integral slope is non-characterizing (\cite[Theorem 1.6(2) and Example 4.5]{bakermotegi}, Figure \ref{fig:bakermotegi}). As a corollary, Theorem \ref{thm:composite} gives the complete list of non-characterizing slopes for these knots.

\begin{cor}
The set of non-characterizing slopes for Baker and Motegi's composite knots consists of all integral slopes.\qed
\end{cor}

To the author's knowledge, this yields the first examples of knots for which the complete list of non-characterizing slopes is known and is not empty. Indeed, the other known examples of a complete list are for the unknot, the trefoil and the figure-8 knot, for which all slopes are characterizing (\cite{kronheimer_monopoles_2007}, \cite{ozsvath_dehn_2019}).

The constraints given by the topology of the exterior of composite knots also lead to the following result.
\begin{theor}\label{thm:cableSFS}
    If $K$ is a knot with an exterior consisting solely of Seifert fibred JSJ pieces, with one of them being a composing space, then any slope that is neither integral nor half-integral is a characterizing slope for $K$.
\end{theor}

\begin{acknow}
I am deeply grateful to my research advisors, Duncan McCoy and Steven Boyer, for their invaluable guidance and the enlightening discussions during which several key ideas were shared. I extend many thanks to Laura Wakelin for several insightful conversations and productive exchanges. I would also like to acknowledge David Futer for bringing to my attention a numerical improvement, and Giacomo Bascapè for his input on the visual aspects of this paper. Lastly, I would like to express my appreciation to the anonymous referee whose comprehensive comments and suggestions greatly improved the quality and clarity of the exposition.
\end{acknow}

\subsection{Structure of paper}
After introducing our notation in Section \ref{sec:notation}, the paper is structured into three main parts. The first, covered in Section \ref{sec:JSJandSurgeredpiece}, describes the JSJ decomposition of a surgery along a knot. The second, consisting of Sections \ref{sec:distslopes}, \ref{sec:sxpiece} and \ref{sec:proof}, presents the proof of Theorem \ref{thm:main}. Finally, Sections \ref{sec:theorexplicitSFS} and \ref{sec:theor2} establish explicit bounds that realize the main theorem for certain families of knots.

\subsection{Outline of proof}
Dehn surgery along a knot $K$ is obtained by gluing a solid torus to the boundary of $S^3_K$, the exterior of $K$ in $S^3$. This boundary is contained in a single JSJ piece of the JSJ decomposition of $S^3_K$. Thus, to understand the topology of a surgery, we must study the fillings of manifolds that arise as JSJ pieces of a knot exterior. We do so in Section \ref{sec:JSJandSurgeredpiece}, where we describe the JSJ decomposition of $S^3_K(p/q)$. In particular, when $|q|$ is sufficiently large, there is one JSJ piece that contains the surgery solid torus; we call it the \textit{surgered piece}.

For a fixed non-trivial knot $K$, suppose there is some knot $K'$ such that there exists an orientation-preserving homeomorphism $S^3_K(p/q) \cong S^3_{K'}(p/q)$. Two scenarios may occur: the surgered piece of $S^3_K(p/q)$ is not mapped to the surgered piece of $S^3_{K'}(p/q)$, or the surgered pieces are mapped one to another. Most of the work towards Theorem \ref{thm:main} lies in the study of the first case. For each possible description of $K$ as a pattern $P$ and a companion knot $J$, we demonstrate that there is a lower bound on $|q|$ determined solely by $K$ such that the outermost JSJ piece of $S^3_J$ is not mapped to the surgered piece of $S^3_{K'}(p/q)$. This yields the following proposition, whose proof occupies Sections \ref{sec:distslopes} and \ref{sec:sxpiece}.

\begin{prop}\label{prop:sxpiece}
Let $K$ be a knot. Suppose $|q|>2$. If there exists an orientation-preserving homeomorphism $S^3_K(p/q) \cong S^3_{K'}(p/q)$ for some knot $K'$, then the homeomorphism sends the surgered piece of $S^3_{K}(p/q)$ to the surgered piece of $S^3_{K'}(p/q)$, provided $|q|$ is sufficiently large. 
\end{prop}

It follows that for $|q|$ sufficiently large, we find ourselves in the situation where an orientation-preserving homeomorphism $S^3_K(p/q) \cong S^3_{K'}(p/q)$ must send the surgered pieces one to another. In that case, the JSJ structures of $S^3_K$ and $S^3_{K'}$ agree away from the JSJ pieces that yielded the surgered pieces. Hence, the problem is now reduced to determining a bound on $|q|$ such that the surgered pieces were in fact obtained from the same manifold. This is done in Section \ref{sec:proof}.

In the final two sections of the paper, we outline explicit bounds that realize Theorem \ref{thm:main} for certain families of knots. 

In Section \ref{sec:theorexplicitSFS}, we provide a lower bound for $|q|$ that ensures that $p/q$ is a characterizing slope for a cable $K$ whose exterior contains only Seifert fibred JSJ pieces. This bound is obtained from the proof of Theorem \ref{thm:main}. In particular, when $K$ is not an $n$-times iterated cable of a torus knot, $n \geq 1$, we show that every slope that is not integral or half-integral is characterizing for $K$.

In Section \ref{sec:theor2}, we demonstrate Theorem \ref{thm:composite}, which gives a realization of Theorem \ref{thm:main} for composite knots when $|q|>1$. Up until this section, we have assumed $|q|>2$, which guaranteed that hyperbolic fillings of JSJ pieces of a knot exterior were also hyperbolic. To lower the bound to $|q|>1$, we need to consider the possibility of exceptional fillings of hyperbolic manifolds. We are able to constrain the topology of half-integer fillings of hyperbolic manifolds of interest by relying on various results that provide upper bounds on the distance between surgery slopes yielding manifolds that contain certain surfaces (\cite{wu_incompressibility_1992}, \cite{gordlue_reducible}, \cite{boyerzhang}), \cite{wu_sutured_1998}). We also use the classification by Gordon and Luecke of hyperbolic knots in $S^3$ and in $S^1 \times D^2$ that admit half-integral toroidal surgeries (\cite{luecke_non-integral_2004}). As a result, we establish that if a non-integral surgery along a knot is obtained from the filling of a hyperbolic JSJ piece, then it can never be orientation-preserving homeomorphic to a non-integral surgery along a composite knot. Theorem \ref{thm:composite} then follows from the argument in Section \ref{sec:proof} regarding composing spaces.

  \section{Notation and preliminaries}\label{sec:notation}
Let $K$ be a knot in $S^3$. We denote by $S^3_K$ the exterior of $K$ in $S^3$, i.e., the manifold obtained by removing the interior of a closed tubular neighbourhood $\nu K$ of $K$ in $S^3$. 
We write $P(J)$ for the satellite knot with pattern $P$ and companion knot $J$. The \textit{winding number} of $P$ is the absolute value of the algebraic intersection number between $P$ and an essential disc in $V = S^1 \times D^2$. The exterior of the satellite $P(J)$ is a gluing $V_P \cup S^3_J$ where $V_P$ denotes the exterior of $P$ seen as a knot in $V$. We call $V_P$ the \textit{pattern space} associated to $P$.

Recall that for any compact irreducible orientable 3-manifold $M$, there is a minimal collection $\mathbf{T}$ of properly embedded disjoint essential tori such that each component of $M \setminus \mathbf{T}$ is either a hyperbolic or a Seifert fibred manifold, and such a collection is unique up to isotopy (\cite{jaco_seifert_1979}, \cite{johannson_homotopy_1979}). The \textit{JSJ decomposition} of $M$ is given by
\[M = M_0 \cup M_1 \cup \ldots \cup M_k,\]
where each $M_i$ is the closure of a component of $M \setminus \mathbf{T}$. A manifold $M_i$ is called a \emph{JSJ piece} of $M$ and a torus in the collection $\mathbf{T}$ is called a \textit{JSJ torus} of $M$. Any homeomorphism between compact irreducible orientable 3-manifolds can be seen as sending JSJ pieces to JSJ pieces, up to isotopy.

The JSJ piece of $S^3_K$ containing the boundary of $\nu K$ is said to be \textit{outermost} in $S^3_K$.

For $\mathcal{T}$ a torus, fix a basis $\{\mu, \lambda\}$ of $H_1(\mathcal{T}; \Z) \cong \Z \oplus \Z$. A simple closed curve on $\mathcal{T}$ represents a class $p\mu+q\lambda$ up to sign, where $p$ and $q$ are coprime.We denote this class by $p/q \in \Q \cup \{1/0\}$ and we call it a \textit{slope}. The \textit{distance} between two slopes $p/q$ and $r/s$ is $\Delta(p/q, r/s) = |ps-qr|$ and it corresponds to the absolute value of the algebraic intersection number between curves representing $p/q$ and $r/s$.

If $M$ is a 3-manifold with toroidal boundary components $\mathcal{T}_1, \ldots, \mathcal{T}_n$ with fixed bases $\{\mu_i, \lambda_i\}$ for each $H_1(\mathcal{T}_i; \Z), i=1, \ldots, n$, then 
\[M(\mathcal{T}_1, \ldots, \mathcal{T}_n ; p_1/q_1, \ldots, p_n/q_n)\]
denotes the Dehn fillings along a simple closed curve representing $p_i/q_i$ on $\mathcal{T}_i$ for each $i=1, \ldots, n$. If only one boundary component of $M$ is filled, we may simply write $M(p/q)$ if it is clear from context which boundary component is filled. If $\del M$ is connected, there is a unique slope $\gamma$ on $\del M$ that has finite order in $H_1(M; \Z)$, called the \textit{rational longitude} of $M$. We refer to the rational longitude as the \textit{longitude} if it is of order 1 in $H_1(M; \Z)$.

When the manifold $M$ is a knot exterior $S^3_K$, a slope $p/q$ on $\del S^3_K$ is expressed in terms of the coordinates of $H_1(\del S^3_K; \Z)$ given by the homotopy class of a curve that bounds an essential disc in $\nu K$, the \textit{meridian of} $S^3_K$, and the homotopy class of a curve that bounds a surface in $S^3_K$, the \textit{longitude of} $S^3_K$, with orientations following the usual convention (a meridional curve pushed into $S^3_K$ and a longitudinal curve have linking number $+1$). The meridian is well-defined by Gordon and Luecke's knot complement theorem (\cite[Theorem 1]{gordlue_knotcompl}) and the longitude is the unique element of $H_1(\del S^3_K; \Z)$ that is null-homologous in $H_1(S^3_K; \Z)$. The slope $1/0$ corresponds to the meridian, while the slope $0/1$ corresponds to the longitude. We will refer to this preferred basis as the one \textit{given by} the knot $K$.

When $K$ is a satellite, we have the following.

\begin{lemma}\label{lemma:pairsPJ}
    Let $K$ be a satellite knot. For each JSJ torus $\mathcal{T}$ of $S^3_K$, there is a pattern $P$ and a knot $J$ such that $K = P(J)$ and $\mathcal{T} = V_P \cap S^3_J$.
\end{lemma}
\begin{proof}
Let $\mathcal{T}$ be a JSJ torus of $S^3_K$. It separates $S^3_K$ into $A\cup_{\mathcal{T}}B$, where $B$ contains $\mathcal{K}=\del S^3_K$. Note that $S^3 \cong S^3_K(1/0) \cong A \cup_{\mathcal{T}} B(\mathcal{K}; 1/0)$. By the loop theorem, any torus in $S^3$ bounds a solid torus, so either $A$ or $B(\mathcal{K}; 1/0)$ must be a solid torus. Since $\mathcal{T}$ is incompressible in $A$ by definition of a JSJ torus, we have that $B(\mathcal{K}; 1/0)$ is a solid torus. Its core is a non-trivial knot $J$ in $S^3$. Thus, $A$ is homeomorphic to $S^3_J$.

    Let $V=B(\mathcal{K}; 1/0)=B \cup_{\mathcal{K}} (\nu K)$. Then $B$ is the solid torus $V$ with the interior of $\nu K$ removed. We can thus see $K$ as a knot in $V$. By the incompressibility of $\mathcal{T}$, $K$ intersects every essential disc in $V$ at least once. Also, $K$ is not the core of $V$ because $\mathcal{T}$ is not boundary parallel in $S^3_K$. Hence, $V \setminus \operatorname{int}(\nu K)$ is the pattern space for a pattern $P$.
\end{proof}

\begin{figure}[h]
    \centering
    \includegraphics[scale=0.4]{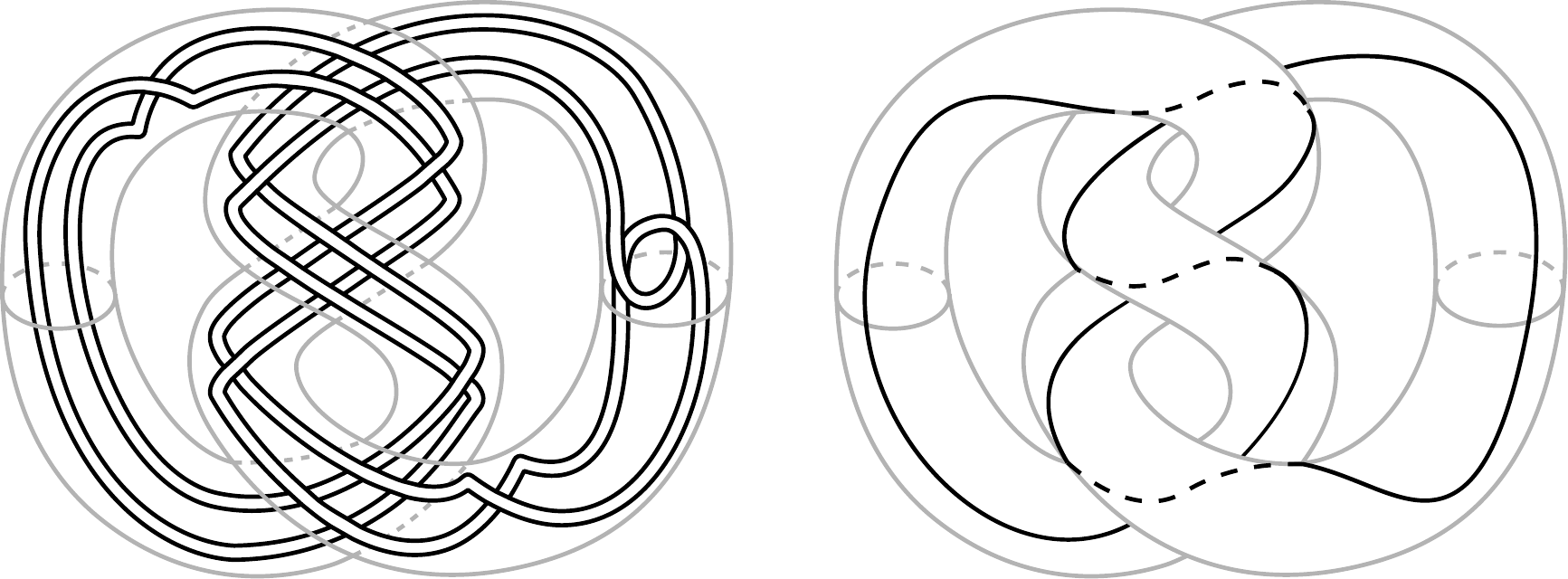}
    \caption{Left: A satellite $K=P(J)$ and the JSJ torus $\mathcal{T} = V_P \cap S^3_J$; Right: The longitude on $\mathcal{T}$ given by $J$}
    \label{fig:longitudeJSJ}
\end{figure}

\begin{defn}\label{def:Tdecomposes}
Let $\mathcal{T}$ be a JSJ torus of $S^3_K$. We say that $\mathcal{T}$ \textit{decomposes} $K$ \textit{into} $P$ \textit{and} $J$ if $\mathcal{T}$ separates $S^3_K$ into $V_P$ and $S^3_J$ as described by Lemma \ref{lemma:pairsPJ}.
\end{defn}
If $\mathcal{T}$ decomposes $K$ into $P$ and $J$, we fix the preferred basis of $H_1(\mathcal{T}; \Z)$ to be the one given by the meridian $\mu_J$ and longitude $\lambda_J$ of $J$ (Figure \ref{fig:longitudeJSJ}), i.e., a slope $p/q$ along $\mathcal{T}$ corresponds to the class $p\mu_J + q \lambda_J \in H_1(\mathcal{T}; \Z) = H_1(\del S^3_J; \Z)$.

Conversely, a pattern space $V_P$ is the data of a knot $P$ in a solid torus $V$, along with a slope $\lambda$ on $\mathcal{T} = \del V$ that intersects $\mu$ once, where $\mu$ is the slope that bounds a disc in $V$. Gluing $V_P$ to a knot exterior $S^3_J$ by respectively identifying $\mu$ and $\lambda$ to the meridian and longitude of $J$ results in the exterior of the knot $K=P(J)$. The preferred basis of $H_1(\mathcal{T}; \Z)$, where $\mathcal{T}$ is seen as a JSJ torus of $S^3_K$ is $\{\mu,\lambda\}$.

Furthermore, for the boundary component $\mathcal{P} = \del \nu P$ of $V_P$, there is a unique class $\lambda_P \in H_1(\mathcal{P}; \Z)$ that is homologous to $w \lambda \in H_1(\mathcal{T}; \Z)$ in $V_P$, where $w$ is the winding number of $P$ (see for instance \cite[p.692]{gord}). The preferred basis of $H_1(\mathcal{P}; \Z)$ is thus given by $\lambda_P$ and $\mu_P$, the class of a curve that bounds an essential disc in $\nu P$.

\begin{lemma}
    Let $K = P(J)$ be a satellite knot. The classes $\lambda_P$ and $\mu_P \in H_1(\mathcal{P}; \Z)$ defined above coincide with the longitude and meridian of $S^3_K$.
\end{lemma}
\begin{proof}
    The meridian of $S^3_K$ and the slope $\mu_P$ coincide because they both bound an essential disc in $\nu K$.

    In $S^3_K$, the class $\lambda_P$ is homologous to $w$ times the longitude $\lambda$ of $S^3_J$, where $w$ is the winding number of $P$. Let $\alpha_P$ be a curve on $\del S^3_K = \mathcal{P}$ representing $\lambda_P$.
    
    If $w=0$, then $\alpha_P$ bounds a surface in $S^3_K$, so $\lambda_P$ is the longitude of $S^3_K$.
    
    If $w\neq0$, then there is a surface $F$ in $S^3_K$ such that 
    \[\del F = (\bigsqcup_{i=1}^w \alpha_i) \sqcup \alpha_P,\]
    where $\alpha_i$ are curves on $\del S^3_J$ representing $\lambda$.
    
    By definition of $\lambda$, each $\alpha_i$ bounds a surface $S_i \subset S^3_J$, $i = 1, \ldots, w$. The union of $F$ with the $S_i$ gives a surface in $S^3_K$ whose boundary is $\alpha_P$. Hence, $\lambda_P$ coincides with the longitude of $S^3_K$.
\end{proof}

If a pattern $P$ in a solid torus $V$ intersects an essential disc in $V$ once, then $P$ is a \textit{composing pattern}. Note that if $K = K_1 \# K_2$ is a composite (or connected sum) of knots $K_1$ and $K_2$, then $K=P_1(K_2)=P_2(K_1)$ where $P_1, P_2$ are composing patterns such that $P_1(U) = K_1$ and $P_2(U) = K_2$, where $U$ is the unknot.

The $(r,s)$-cable of a knot $J$ is denoted by $C_{r,s}(J)$, where $s$ is the winding number of the cable pattern. We may assume that $s>0$ since the $(r,s)$ and $(-r,-s)$-cable patterns are equivalent.
The pattern space $V_{C_{r,s}}$ is an \textit{$(r,s)$-cable space}. It is the outermost JSJ piece of the exterior of $C_{r,s}(J)$. Further, it admits a Seifert fibration with base orbifold an annulus with one cone point of order $s$. On its boundary component corresponding to $\del S^3_{C_{r,s}(J)}$, the $(r,s)$-cable space has regular fibres of slope $rs/1$. On the other boundary component coinciding with $\del S^3_{J}$, a regular fibre has slope $r/s$ in the coordinates given by $J$.

We denote by $T_{a,b}$ the $(a,b)$-torus knot. Its exterior is Seifert fibred, with two exceptional fibres of orders $|a|$ and $|b|$. The regular fibres have slope $ab/1$ on $\del S^3_{T_{a,b}}$.

  \section{JSJ decompositions and the surgered piece}\label{sec:JSJandSurgeredpiece}

The JSJ pieces of a non-trivial knot exterior take on one of four special types. Here is a version of this result found in \cite{bud}.

\begin{thm}[{\cite[Theorem 4.18]{bud}}]\label{thm:budney}
In the JSJ decomposition of the exterior of a non-trivial knot, the outermost JSJ piece is either
\begin{enumerate}
    \item the exterior of a torus knot;
    \item a composing space, i.e., a Seifert fibre space with at least 3 boundary components and base orbifold a planar surface with no cone points;
    \item the exterior of a hyperbolic knot or link such that if the component of the link corresponding to the knot is removed, the resulting link is the unlink;
    \item a cable space, i.e., a Seifert fibre space with base orbifold an annulus with one cone point.
\end{enumerate}
\end{thm}

By Lemma \ref{lemma:pairsPJ}, a JSJ torus of the exterior $S^3_K$ of a knot $K$ is the boundary of the exterior of a non-trivial knot in $S^3$. Therefore, each JSJ piece of $S^3_K$ is the outermost piece of some knot exterior, which implies that every JSJ piece of $S^3_K$ belong to one of the types listed in Theorem \ref{thm:budney}. 

Homological calculations from \cite{gord} lead to the following two results.

\begin{lemma}[{\cite[Lemma 3.3]{gord}}]\label{lemma:nullhom}
Let $P(J)$ be a satellite knot, where $P$ has winding number $w$. Denote the boundary components of the pattern space $V_P$ by $\mathcal{P} = \del \nu P$ and $\mathcal{T}=\del S^3_J$.
    \begin{enumerate}[(i)]
        \item $H_1(V_P(\mathcal{P};p/q); \Z) \cong \Z \oplus (\Z/g_{p,w}\Z)$, where $g_{p,w}$ is the greatest common divisor of $p$ and $w$;
        \item The kernel of $H_1(\mathcal{T}; \Z) \rightarrow H_1(V_P(\mathcal{P};p/q); \Z)$ induced by inclusion is generated by
        \[\begin{cases}
            (p/g_{p,w}) \mu + (qw^2/g_{p,w}) \lambda & \text{if }w\neq0 \\
            \mu     & \text{if }w=0
        \end{cases},\]
        where $\{\mu, \lambda\}$ is the basis of $H_1(\mathcal{T}; \Z)$ given by $J$.
    \end{enumerate}
\end{lemma}

\begin{prop}[{\cite[Corollary 7.3]{gord}}]\label{prop:cablepqs2}
    Let $K = C_{r,s}(J)$ be cable knot. 
    \begin{enumerate}
    \item If $|qrs-p|>1$, then $S^3_{K}(p/q)$ is the union along their boundary of $S^3_J$ and a Seifert fibre space with incompressible boundary; 
    \item If $|qrs-p|=1$, then $S^3_{K}(p/q) \cong S^3_J(p/(qs^2))$.
    \end{enumerate}
\end{prop}
Note that if $|qrs-p|=1$, then $g_{p,s}=1$ and $p/(qs^2)$ is a well-defined slope.

Gordon and Luecke showed that if $p/q$ is not an integer, then the surgery $S^3_K(p/q)$ is irreducible (\cite[Theorem 1]{gorlue_integral_1987}). Thus, it admits a JSJ decomposition. For the rest of this section, we focus our attention on the topology of the JSJ pieces of $S^3_K(p/q)$ when $|q|>2$. The next theorem from \cite{lack} combines results from various authors (\cite{gordon_toroidal_1999}, \cite{gordon_annular_2000}, \cite{scharlemann_producing_1990}), \cite{wu_incompressibility_1992}).

\begin{thm}[{\cite[Theorem 2.8]{lack}}]\label{thm:hyplink}
Let $M$ be the exterior of a hyperbolic link in $S^3$ with components $L_0, L_1, \ldots, L_n, n \geq 1$, such that the link formed by the components $L_1, \ldots, L_n$ is the unlink. Let $\sigma$ be a slope on $\mathcal{L}_0 = \del\nu L_0 \subset \del M$ and let $\mu$ be the slope on $\mathcal{L}_0$ that bounds a disc in $\nu L_0$. If $\Delta(\sigma, \mu) > 2$, then $M(\mathcal{L}_0; \sigma)$ is hyperbolic.
\end{thm}

Let $K$ be a non-trivial knot and $Y_0 \cup Y_1 \cup \ldots \cup Y_k$ be the JSJ decomposition of its exterior $S^3_K$, where $Y_0$ is the outermost piece. The Dehn surgery $S^3_K(p/q)$ is obtained by filling $Y_0$ along $\mathcal{K} = \del S^3_K \subset \del Y_0$.

\begin{prop}\label{prop:fillingY0}
    If $|q|>2$, the filling $Y_0(\mathcal{K}; p/q)$ is either a Seifert fibre space or a hyperbolic manifold. In particular,
    \begin{enumerate}
        \item If $Y_0$ is the exterior of a hyperbolic link that is not a knot, then $Y_0(\mathcal{K}; p/q)$ is hyperbolic;
        \item If $Y_0$ is a composing space, then $Y_0(\mathcal{K}; p/q)$ is Seifert fibred with base orbifold a planar surface with at least two boundary components and one cone point of order $|q|$;
        \item If $Y_0$ is an $(r,s)$-cable space and $|qrs-p|>1$, then $Y_0(\mathcal{K}; p/q)$ is Seifert fibred with base orbifold a disc with two cone points of orders $|qrs-p|$ and $s$;
        \item If $Y_0$ is an $(r,s)$-cable space and $|qrs-p|=1$, then $Y_0(\mathcal{K}; p/q)$ is a solid torus.
    \end{enumerate}
\end{prop}
\begin{proof}
If $Y_0 = S^3_K$ and $K$ is a hyperbolic knot, then $Y_0(\mathcal{K}; p/q) = S^3_K(p/q)$ does not contain an essential sphere or an incompressible torus if $|q|>2$, so it is either hyperbolic or Seifert fibred (\cite[Theorem 1.1]{gordon_dehn_1995}).
If $Y_0 = S^3_K$ and $K$ is a torus knot, then $Y_0(\mathcal{K}; p/q) = S^3_K(p/q)$ is Seifert fibred if $|q|>1$ (\cite[Proposition 3.1]{moser}).

If $Y_0$ is the exterior of a hyperbolic link, then by Theorem \ref{thm:hyplink}, $Y_0(\mathcal{K}; p/q)$ is hyperbolic if $|q|>2$.

If $Y_0$ is a composing space, a regular fibre on $\mathcal{K}$ has slope $1/0$. If $|q| > 1$, we have $\Delta(1/0, p/q) = |q|>1$, so the surgery slope does not coincide with the regular fibre slope. The Seifert fibred structure of $Y_0$ thus extends to the surgery solid torus adding an exceptional fibre of order $|q|$. Moreover, a composing space has at least three boundary components, so $Y_0(\mathcal{K}; p/q)$ has at least two boundary components.

If $Y_0$ is an $(r,s)$-cable space, a regular fibre on $\mathcal{K}$ has slope $rs/1$. If $|q| > 1$, we have $\Delta(rs/1,p/q) = |qrs-p|\neq 0$, so the surgery slope does not coincide with the regular fibre slope. The Seifert fibred structure of $Y_0$ thus extends to the surgery solid torus. If $|qrs-p|>1$, the surgery adds an exceptional fibre of order $|qrs-p|$. If $|qrs-p|=1$, then the surgery solid torus is regularly fibred in $Y_0(\mathcal{K}; p/q)$, so $Y_0(\mathcal{K}; p/q)$ has base orbifold a disc and one cone point. It is a solid torus.
\end{proof}

\begin{prop}\label{prop:JSJsurgery}
Suppose $|q|>2$. The JSJ decomposition of $S^3_K(p/q)$ is either
\[Y_0(\mathcal{K}; p/q) \cup Y_1 \cup Y_2 \cup \ldots \cup Y_k\]
or
\[Y_1(\mathcal{J}; p/(qs^2)) \cup Y_2 \cup \ldots \cup Y_k,\]
where $\mathcal{J} = Y_0 \cap Y_1$ and $s\geq 2$. The second scenario occurs precisely when $K$ is a cable knot $C_{r,s}(J)$, $Y_1$ is the outermost piece of $S^3_J$, and $|qrs-p|=1$.
\end{prop}
\begin{proof}
    By the previous proposition, $Y_0(\mathcal{K}; p/q)$ is either Seifert fibred or hyperbolic. If it is hyperbolic or closed, then the result is immediate.

    If $Y_0(\mathcal{K}; p/q)$ is Seifert fibred and has boundary, i.e., in cases (2),(3) and (4) of Proposition \ref{prop:fillingY0}, then $Y_0(\mathcal{K}; p/q)$ might admit a Seifert structure that extends across adjacent JSJ pieces. By definition of the JSJ decomposition, this structure would have to differ from the one inherited from the Seifert structure on $Y_0$.

    Only cases (3) and (4) of Proposition \ref{prop:fillingY0}, which correspond to $K$ being the cable of a knot $J$,   
    may give rise to manifolds $Y_0(\mathcal{K}; p/q)$ that admit multiple Seifert fibred structures.
    
    In case (3), $Y_0(\mathcal{K}; p/q)$ admits more than one Seifert fibred structure when it is a twisted $I$-bundle over the Klein bottle. One is inherited from $Y_0$ and has base orbifold a disc with two cone points each of order 2, and the other has base orbifold a Möbius band with no cone points. 
    
    This second structure has regular fibres that are non-meridional and non-integral on $\del Y_0(\mathcal{K}; p/q)$ if $|q|>1$. Indeed, a regular fibre of this structure corresponds to the generator of the $\Z/2\Z$ summand of $H_1(Y_0(\mathcal{K};p/q); \Z) \cong \Z \oplus \Z/2\Z$. Let $p'\mu_J+q'\lambda_J \in  H_1(\del Y_0(\mathcal{K};p/q); \Z)$ be the class of a regular fibre in the coordinates given by $J$. Let $i : H_1(\del Y_0(\mathcal{K};p/q); \Z) \rightarrow H_1(Y_0(\mathcal{K};p/q); \Z)$ be induced by inclusion. Then $2(p'\mu_J+q'\lambda_J)$ generates the kernel of $i$. Since $s=|qrs-p| = 2$, $p$ is even. By Lemma \ref{lemma:nullhom}(ii), $\ker i$ is generated by $(p/2)\mu_J + 2q\lambda_J$. Since $q$ and $r$ are odd, we have that $p$ is divisible by 4, from which we deduce that $p'\mu_J+q'\lambda_J = (p/4)\mu_J + q\lambda_J$, which is non-meridional and non-integral.

    Therefore, this Seifert fibred structure does not extend to an adjacent Seifert fibred JSJ piece $Y_1$, because the slope of a regular fibre of $Y_1$ on the JSJ torus $\mathcal{J} = Y_0 \cap Y_1$ is either meridional (if $Y_1$ is a composing space) or integral (if $Y_1$ is a torus knot exterior or a cable space) in the coordinates given by the meridian $\mu_J$ and longitude $\lambda_J$ of the companion knot $J$.

    In case (4), $K$ is a cable knot $C_{r,s}(J)$ such that $|qrs-p|=1$, and $Y_0(\mathcal{K}; p/q)$ is a solid torus. By Proposition \ref{prop:cablepqs2}, $S^3_K(p/q) \cong S^3_J(p/(qs^2))$. We have $|qs^2| > |q| > 2$. We iterate the above argument for $S^3_J(p/(qs^2))$ to reduce to case (4) of Proposition \ref{prop:fillingY0} for $S^3_J(p/(qs^2))$. We show that this case does not occur if $|q|>1$.

    Suppose $Y_1(\mathcal{J}; p/(qs^2))$ is a solid torus. Then $Y_1$ must be an $(r', s')$-cable space and $|qrs-p|=|qs^2r's' - p| = 1$ (Proposition \ref{prop:cablepqs2}). Hence, $|q(rs - s^2r's')| = 2$ or $0$. As $|q|, s > 1$, the first case does not occur, and the second case happens only if $rs - s^2r's' = 0$, but this contradicts $r$ and $s$ being coprime.
\end{proof}

It follows that the surgery solid torus is contained in exactly one JSJ piece of $S^3_K(p/q)$ when $|q|>2$.

\begin{defn}\label{def:sxpiece}
    Suppose $|q|>2$. The \emph{surgered piece} of $S^3_K(p/q)$ is the JSJ piece of $S^3_K(p/q)$ that contains the surgery solid torus. It corresponds to either $Y_0(\mathcal{K}; p/q)$ or $Y_1(\mathcal{J}; p/(qs^2))$, as outlined in Proposition \ref{prop:JSJsurgery}.
\end{defn}

\color{black}


The topology of the surgered piece is summarized as follows.

\begin{prop}\label{prop:geomsxpiece}
    Suppose $|q|>2$. The surgered piece of $S^3_K(p/q)$ is a filling $Y(p/(qt^2))$ of a JSJ piece $Y$ of $S^3_K$, for some integer $t \geq 1$. In particular,
    \begin{enumerate}
        \item $Y(p/(qt^2))$ has non-empty boundary and is hyperbolic if and only if $Y$ is the exterior of a hyperbolic link that is not a knot;
        \item $Y(p/(qt^2))$ is Seifert fibred with base orbifold a planar surface with at least two boundary components and one cone point of order $|qt^2|$ if and only if $Y$ is a composing space;
        \item$Y(p/(qt^2))$ is Seifert fibred with base orbifold a disc with two cone points if and only if $Y$ is a cable space. In particular, if $Y$ is an $(r,s)$-cable space, then the cone points have orders $|qt^2rs-p|>1$ and $s$.
    \end{enumerate}
    Furthermore, if $|q|>8$, then
    \begin{enumerate}\setcounter{enumi}{3}
        \item $Y(p/(qt^2))$ is closed and Seifert fibred if and only if $Y$ is the exterior of a torus knot;
        \item $Y(p/(qt^2))$ is closed and hyperbolic if and only if $Y$ is the exterior of a hyperbolic knot.
    \end{enumerate}
\end{prop}
\begin{proof}
    The converses of (1), (2), (3) follow from Proposition \ref{prop:fillingY0}. We deduce the direct implications from Theorem \ref{thm:budney} as follows. 
    
    If $Y(p/(qt^2))$ is not closed and is hyperbolic, then $Y$ is hyperbolic with at least two boundary components, so it must be the exterior of a hyperbolic link that is not a knot. 
    
    If $Y(p/(qt^2))$ is Seifert fibred and has $n \geq 1$ boundary components, then $Y$ must be Seifert fibred (Theorem \ref{thm:hyplink}) and it has $n+1$ boundary components. Hence, if $n \geq 2$, $Y$ is a composing space, while if $n=1$, $Y$ is a cable space.
    
    When $|q|>8$, it is a result of Lackenby and Meyerhoff (\cite[Theorem 1.2]{lackmey}) that if $Y$ is the exterior of a hyperbolic knot, then $Y(p/(qt^2))$ must also be hyperbolic. Conversely, if $Y(p/(qt^2))$ is closed and hyperbolic, then $Y$ is the exterior of knot that must be hyperbolic.
    
    If $Y$ is the exterior of a torus knot, then $Y(p/(qt^2))$ is Seifert fibred (\cite[Proposition 3.1]{moser}). Conversely, if $Y(p/(qt^2))$ is closed and Seifert fibred, $Y$ is a knot exterior that must be Seifert fibred, by the result of Lackenby and Meyerhoff. The only knots whose exteriors are Seifert fibred are torus knots (\cite[Theorem 2]{moser}).    
\end{proof}

The five types of surgered pieces described in Proposition \ref{prop:geomsxpiece} correspond to fillings of distinct types of JSJ pieces of a knot exterior.
\begin{cor}\label{cor:sxtype}
    Suppose $|q|>2$ and let $K$ and $K'$ be knots. Suppose further that the surgered piece $Y(p/(qt^2))$ of $S^3_K(p/q)$ is homeomorphic to the surgered piece $Y'(p'/(q'(t')^2))$ of $S^3_{K'}(p'/q')$. 
    \begin{enumerate}
        \item     If $Y(p/(qt^2))$ and $Y'(p'/(q'(t')^2))$ have non-empty boundary, then $Y$ and $Y'$ are of the same type, as listed by Theorem \ref{thm:budney}.
        \item     Furthermore, if $|q|>8$ and if $Y(p/(qt^2))$ and $Y'(p'/(q'(t')^2))$ are closed, then $Y$ and $Y'$ are both torus knots or both hyperbolic knots.\qed
    \end{enumerate}
\end{cor}

Comparing with Theorem \ref{thm:budney}, we obtain additional constraints on the structure of the surgered piece.

\begin{prop}\label{prop:geomswitch}
 Suppose $|q|>2$. Let $Y$ be the JSJ piece of $S^3_K$ such that the surgered piece of $S^3_K(p/q)$ is a filling $Y(p/(qt^2))$ for some integer $t \geq 1$. If $Y(p/(qt^2))$ is homeomorphic to a JSJ piece of a knot exterior, then
 \begin{enumerate}
     \item $Y$ is not the exterior of a knot;
     \item $Y(p/(qt^2))$ is homeomorphic to the exterior of a hyperbolic knot or link such that if a specific component of the link is removed, the resulting link is the unlink, if and only if $Y$ is hyperbolic;
     \item $Y(p/(qt^2))$ is homeomorphic to an $(r, |qt^2|)$-cable space if and only if $Y$ is a composing space;
     \item $Y(p/(qt^2))$ is homeomorphic to the exterior of a torus knot if and only if $Y$ is an $(r,s)$-cable space.
 \end{enumerate}
\end{prop}
\begin{proof}
    For (1), we observe that if $Y$ is the exterior of a knot, then $Y(p/(qt^2))$ is a closed manifold. However, all JSJ pieces of a knot exterior have non-empty boundary.

    The implications of (2), (3) and (4) follow from Proposition \ref{prop:geomsxpiece}. We show their converses.

    If $Y$ is hyperbolic, then by (1), it is not the exterior of a knot. Hence, $Y(p/(qt^2))$ is hyperbolic by Proposition \ref{prop:geomsxpiece}. By Theorem \ref{thm:budney}, a hyperbolic JSJ piece of a knot exterior is as stated in (2).

    If $Y$ is a composing space, then $Y(p/(qt^2))$ is Seifert fibred with only one exceptional fibre of order $|qt^2|$ (Proposition \ref{prop:geomsxpiece}). By Theorem \ref{thm:budney}, cable spaces are the only Seifert fibred JSJ pieces of a knot exterior with only one exceptional fibre. An $(r,s)$-cable space has an exceptional fibre of order $s$, so $Y(p/(qt^2))$ is an $(r, |qt^2|)$-cable space.
    
    If $Y$ is an $(r,s)$-cable space, $Y(p/(qt^2))$ is Seifert fibred with two exceptional fibres (Proposition \ref{prop:geomsxpiece}). By Theorem \ref{thm:budney}, torus knot exteriors are the only Seifert fibred JSJ pieces of a knot exterior with two exceptional fibres.
\end{proof}

  \section{Distinguished slopes}\label{sec:distslopes}
The goal of Sections \ref{sec:distslopes} and \ref{sec:sxpiece} is to prove the following proposition.

\begin{prop}\label{prop:sxpiecesbound}
Let $K$ be a satellite knot and $\mathcal{T}$ be a JSJ torus of $S^3_K$ that decomposes $K$ into $P$ and $J$. There exists a constant $L(\mathcal{T})$ with the following property. Suppose $|q|>2$. If there exists an orientation-preserving homeomorphism $S^3_K(p/q) \cong S^3_{K'}(p/q)$ for some knot $K'$, then the homeomorphism does not map the outermost piece of $S^3_{J} \subset S^3_K(p/q)$ to the surgered piece of $S^3_{K'}(p/q)$, provided $|q| > L(\mathcal{T})$.
\end{prop} 

Note that if $\mathcal{T}$ is compressible in the statement of Proposition \ref{prop:sxpiecesbound}, then the outermost piece of $S^3_J$ inside $S^3_K(p/q)$ is not a JSJ piece of $S^3_K(p/q)$, so the conclusion holds with $L(\mathcal{T})=2$. Therefore, in the subsequent discussion, we will assume that $\mathcal{T}$ is incompressible in $S^3_K(p/q)$.

\subsection{Filled patterns and companion knots}\label{sec:scenario}
Throughout Section \ref{sec:distslopes}, we will consider the following scenario.

The satellite knot $K$ is fixed. Suppose there is a knot $K'$ such that there exists an orientation-preserving homeomorphism $S^3_K(p/q) \cong S^3_{K'}(p/q)$, mapping JSJ pieces of $S^3_K(p/q)$ to JSJ pieces of $S^3_{K'}(p/q)$.

Let $X$, $X'$ be the surgered pieces of $S^3_{K}(p/q)$, $S^3_{K'}(p/q)$ respectively. Suppose that the homeomorphism does not send $X$ to $X'$. Then $X'$ is the image of a JSJ piece of $S^3_{K}(p/q)$ that is not $X$. That JSJ piece in $S^3_{K}(p/q)$ is the outermost piece of $S^3_J \subset S^3_K$ for some knot $J$. Let $\mathcal{T} = \del S^3_{J}$. This is a JSJ torus of $S^3_{K}$. By Lemma \ref{lemma:pairsPJ}, $\mathcal{T}$ decomposes $K$ into a pattern $P$ and the knot $J$.

The JSJ torus $\mathcal{T}$ is sent by the homeomorphism to a JSJ torus $\mathcal{T}'$ of $S^3_{K'}(p/q)$, which is also a JSJ torus of $S^3_{K'}$ by Proposition \ref{prop:JSJsurgery}. By Lemma \ref{lemma:pairsPJ}, $\mathcal{T}'$ decomposes $K'$ into a pattern $P'$ and a knot $J'$.

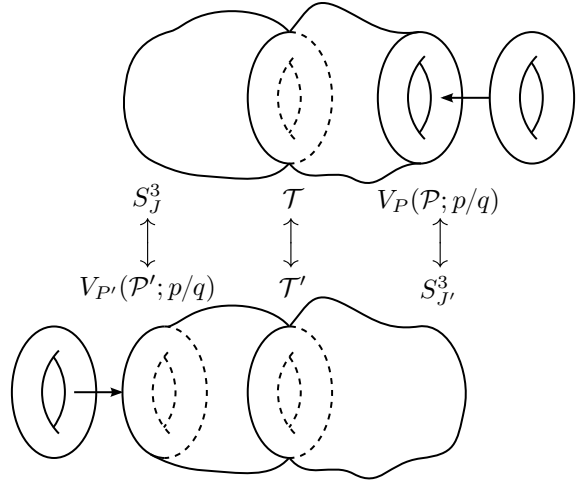
\begin{figure}[h]
    \centering
    \begin{tikzpicture}
    \node at (0,0) {\includegraphics[scale=0.5]{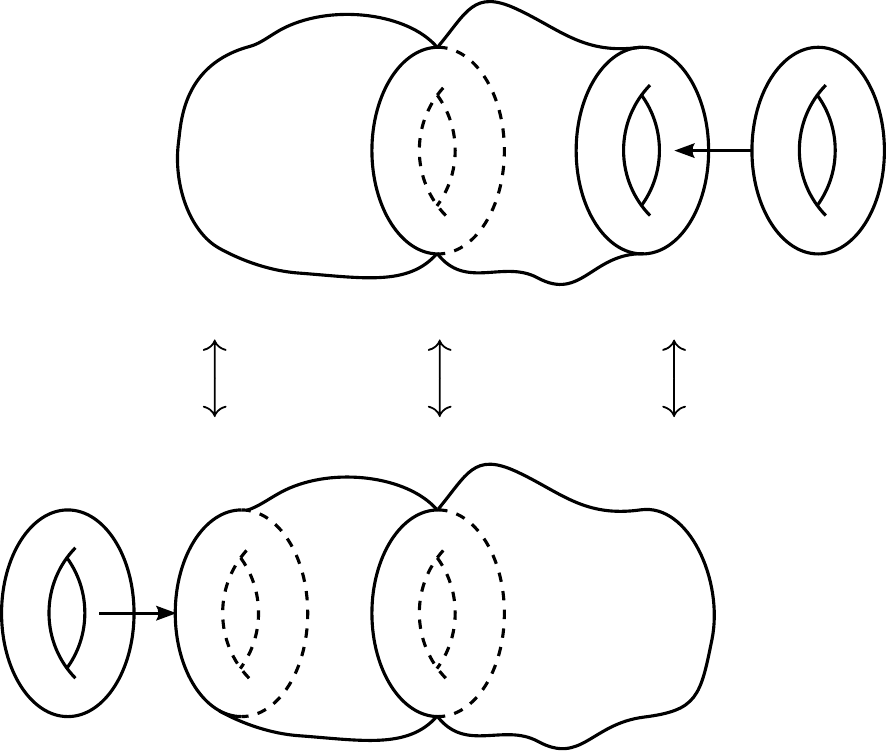}};

    \node at (0,16pt) {$\mathcal{T}$};
    \node at (55pt,15pt) {$V_{P}(\mathcal{P};p/q)$};
    \node at (-55pt,16pt) {$S^3_{J}$};   
    
    \node at (0,-17pt) {$\mathcal{T}'$};
    \node at (-55pt,-18pt) {$V_{P'}(\mathcal{P}';p/q)$};
    \node at (55pt,-19pt) {$S^3_{J'}$};
    
    \end{tikzpicture}
    \caption{Homeomorphism carrying filled pattern space to exterior of companion knot}
    \label{fig:switchinghomeo}
\end{figure}

 Let $\mathcal{P}$ and $\mathcal{P}'$ respectively denote the boundary components $\del \nu P$ and $\del \nu P'$ of $V_P$ and $V_{P'}$, the pattern spaces associated to $P$ and $P'$. The homeomorphism $S^3_K(p/q) \cong S^3_{K'}(p/q)$ then restricts to a homeomorphism between $V_P(\mathcal{P};p/q)$ and $S^3_{J'}$, and between $S^3_J$ and $V_{P'}(\mathcal{P}';p/q)$ (Figure \ref{fig:switchinghomeo}).

We will now identify distinguished slopes on the JSJ tori $\mathcal{T} \subset S^3_K(p/q)$ and $\mathcal{T}' \subset S^3_{K'}(p/q)$. Information about the gluing of JSJ pieces along their boundaries will be obtained by analyzing the distances between these slopes. In Section \ref{sec:sxpiece}, we will rely on the fact that distances between slopes are preserved by homeomorphisms to establish constraints on the coefficients $p$ and $q$.

\subsection{General pattern case}
In the scenario described in Section \ref{sec:scenario} and Figure \ref{fig:switchinghomeo}, we have the following lemma.

\begin{lemma}\label{lemma:gw1}
The greatest common divisor $g_{p,w}$ of $p$ and $w$ is 1.
\end{lemma}
\begin{proof}
On one hand, we have $H_1(V_P(\mathcal{P};p/q);\Z) = \Z \oplus (\Z/g_{p,w}\Z)$ (Lemma \ref{lemma:nullhom}(i)). On the other hand, $H_1(S^3_{J'};\Z) = \Z$. Since $S^3_{J'} \cong V_P(\mathcal{P};p/q)$, we conclude that $g_{p,w} = 1$.
\end{proof}

Our first distinguished slope on $\mathcal{T} \subset S^3_K(p/q)$ is the longitude of $V_P(\mathcal{P};p/q)$ seen as a knot exterior. Combining Lemma \ref{lemma:nullhom}(ii) with Lemma \ref{lemma:gw1}, we have that this slope corresponds to the class
        \[\begin{cases}
            p\mu_J + qw^2\lambda_J & \text{if }w\neq0, \\
            \mu_J     & \text{if }w=0
        \end{cases}\]
in $H_1(\mathcal{T}; \Z) = H_1(\del S^3_J; \Z)$, where $\mu_J$ and $\lambda_J$ are the meridian and longitude of $S^3_J$ respectively.

Our second distinguished slope on $\mathcal{T} \subset S^3_K(p/q)$ is the meridian of $V_P(\mathcal{P}; p/q)$ seen as a knot exterior. Let $x\mu_J + y\lambda_J$ be a class corresponding to this slope.

On $\mathcal{T}'\subset S^3_{K'}(p/q)$, we have two analogous distinguished slopes: the meridian and the longitude of $V_{P'}(\mathcal{P}';p/q)$ seen as a knot exterior. Let $x'\mu_{J'} + y'\lambda_{J'}$ be a class in $H_1(\mathcal{T}'; \Z) = H_1(\del S^3_{J'}; \Z)$ corresponding to this meridian, where $\mu_{J'}$ and $\lambda_{J'}$ are the meridian and longitude of $S^3_{J'}$ respectively.

\begin{lemma}\label{lemma:meridian}
\begin{enumerate}[(i)]
    \item The meridian $x\mu_J + y\lambda_J$ of $V_P(\mathcal{P};p/q)$ is such that $|x| =  |q(w')^2|$.
    \item The meridian $x'\mu_{J'} + y'\lambda_{J'}$ of $V_{P'}(\mathcal{P}';p/q)$ is such that $|x'| =  |qw^2|$.
\end{enumerate} 
\end{lemma}
\begin{proof}
The homeomorphism $S^3_K(p/q) \cong S^3_{K'}(p/q)$ sends the meridian $x\mu_J + y\lambda_J$ of $V_P(\mathcal{P};p/q)$ to the meridian $\mu_{J'}$ of $S^3_{J'}$, and the longitude $\lambda_J$ of $S^3_J$ to the longitude $p\mu_{J'}+q(w')^2\lambda_{J'}$ of $V_{P'}(\mathcal{P}';p/q)$. Since the homeomorphism preserves distances between slopes, we have
\[|x| = \Delta\left(\frac{x}{y}, \frac{0}{1}\right) = \Delta\left(\frac{1}{0}, \frac{p}{q(w')^2}\right) = |q(w')^2|.\]
We obtain (ii) symmetrically. 
\end{proof}

Table \ref{table:distslopes} summarizes this discussion.

\begin{table}[h]
    \centering
    \begin{tabular}{|c|c|c|}
    \hline
    & $\mathcal{T}$ & $\mathcal{T}'$ \\ \hline\hline
       \begin{tabular}{c}
     Meridian of  \\
     $V_P(\mathcal{P};p/q) \cong S^3_{J'}$
 \end{tabular}  & 
 \multicolumn{1}{l|}{
    \begin{tabular}{@{}l@{}l}
        $q(w')^{2}\mu _{J}+y\lambda _{J}$ & $\text{ if } w' \neq 0$ \\
        $\lambda _{J}$ & $\text{ if } w'=0$
    \end{tabular}
    }
    & $\mu_{J'}$\\\hline
     \begin{tabular}{c}
     Longitude of  \\
     $V_P(\mathcal{P};p/q) \cong S^3_{J'}$
 \end{tabular}  &  
 \multicolumn{1}{l|}{
    \begin{tabular}{@{}l@{}l}
        \rlap{$p\mu _{J}+qw^{2}\lambda _{J}$}\phantom{$q(w')^{2}\mu _{J}+y\lambda _{J}$} & $\text{ if } w \neq 0$ \\
        $\mu _{J}$ & $\text{ if } w=0$
    \end{tabular}
    }
    & $\lambda_{J'}$ \\\hline
 \begin{tabular}{c}
      Meridian of  \\
     $S^3_J \cong V_{P'}(\mathcal{P}';p/q)$
 \end{tabular} & $\mu_J$ & 
 \multicolumn{1}{l|}{
     \begin{tabular}{@{}l@{}l}
	\rlap{$qw^{2}\mu _{J'}+y'\lambda _{J'}$}\phantom{$p\mu _{J'}+q(w')^{2}\lambda _{J'}$} &$\text{ if } w \neq 0 $\\
    $\lambda _{J'}$ &$\text{ if } w=0$
\end{tabular}
}
 \\\hline
 \begin{tabular}{c}
     Longitude of  \\
     $S^3_J \cong V_{P'}(\mathcal{P}';p/q)$
 \end{tabular} & $\lambda_J$ & 
 \multicolumn{1}{l|}{
    \begin{tabular}{@{}l@{}l}
    $p\mu _{J'}+q(w')^{2}\lambda _{J'}$ & $\text{ if } w' \neq 0$ \\
    $\mu _{J'}$ & $\text{ if } w'=0$
\end{tabular}
}
 \\ \hline
    \end{tabular}
    \caption{Distinguished slopes on $\mathcal{T}$ and $\mathcal{T}'$}
    \label{table:distslopes}
\end{table}

\subsection{Iterated cable case}\label{sec:distslopesCT}

In the case where $P$ is an iterated cable, we also distinguish the slopes of regular Seifert fibres in the scenario described in Section \ref{sec:scenario} and Figure \ref{fig:switchinghomeo}. 

\begin{lemma}
    If $P$ is an iterated cable $C_{r_n, s_n} \ldots C_{r_2, s_2}C_{r_1, s_1}, n\neq 1$, then
    \begin{enumerate}
        \item The JSJ piece of $V_P(\mathcal{P}; p/q)$ with boundary component $\mathcal{T}$ is Seifert fibred and its regular fibre has slope $r_1\mu_J + s_1\lambda_J$ on $\mathcal{T}$;
        \item The outermost JSJ piece of $S^3_{J'}$ is Seifert fibred and its regular fibre has integral slope $k\mu_{J'} + \lambda_{J'}$ on $\mathcal{T}'$, for some $k \in \Z$.
    \end{enumerate}

\end{lemma}
\begin{proof}
Let $V_i$ be $(r_i, s_i)$-cable spaces for $i = 1, \ldots, n$. The pattern space $V_P$ has JSJ decomposition $V_1 \cup V_2 \cup \ldots \cup V_n$, where $\mathcal{T} \subset \del V_1$. A regular fibre of $V_1$ has slope $r_1\mu_J + s_1\lambda_J$ on $\mathcal{T}$. 

If $V_P(\mathcal{P}; p/q)$ contains an incompressible torus, then it is clear that the regular fibre slope on $\mathcal{T}$ remains unchanged. Furthermore, $V_1$ is homeomorphic to the outermost piece of $S^3_{J'}$, which must also be a cable space. Hence, a regular fibre of the outermost piece of $S^3_{J'}$ has integral slope on $\mathcal{T}' = \del S^3_{J'}$.

If $V_P(\mathcal{P}; p/q)$ contains no incompressible torus, then it is a filling of $V_1$ by Proposition \ref{prop:cablepqs2}. By hypothesis (Figure \ref{fig:switchinghomeo}), this filling is homeomorphic to a JSJ piece of $S^3_{K'}$. By Proposition \ref{prop:geomswitch}, this piece is the exterior of a torus knot. It follows that the Seifert fibred structure on $V_P(\mathcal{P}; p/q)$ is unique, and it is the one inherited from $V_1$. Moreover, the torus $\mathcal{T}' \subset S^3_{K'}(p/q)$ is the boundary of a torus knot exterior, so a regular fibre has integral slope on $\mathcal{T}'$.
\end{proof}

Thus, the homeomorphism $S^3_K(p/q) \cong S^3_{K'}(p/q)$ maps the slope $r_1\mu_J + s_1\lambda_J$ on $\mathcal{T}$ to a slope $k\mu_{J'} + \lambda_{J'}$ on $\mathcal{T}'$, where $k \in \Z$.

\begin{table}[h]
    \centering
    \begin{tabular}{|c|c|c|}
    \hline
    & $\mathcal{T}$ & $\mathcal{T}'$ \\ \hline\hline
 \begin{tabular}{c}
     Regular fibre in  \\
     $V_P(\mathcal{P}; p/q) \cong S^3_{J'}$
 \end{tabular} & $r_1\mu_J + s_1\lambda_J$ & $k\mu_{J'} + \lambda_{J'}$  \\
 \hline
    \end{tabular}
    \caption{Regular fibre slopes on $\mathcal{T}$ and $\mathcal{T}'$}
    \label{table:distslopesCT}
\end{table}

  \section{Surgered pieces are sent to surgered pieces}\label{sec:sxpiece}

This section is dedicated to demonstrating Proposition \ref{prop:sxpiecesbound}, from which Proposition \ref{prop:sxpiece} follows easily.

\begin{prop:sxpiece}
Let $K$ be a knot. Suppose $|q|>2$. If there exists an orientation-preserving homeomorphism $S^3_K(p/q) \cong S^3_{K'}(p/q)$ for some knot $K'$, then the homeomorphism sends the surgered piece of $S^3_{K}(p/q)$ to the surgered piece of $S^3_{K'}(p/q)$, provided $|q|$ is sufficiently large. 
\end{prop:sxpiece}

	\begin{proof} 
If $K$ is not a satellite knot, then the result follows from Proposition \ref{prop:JSJsurgery}, so suppose $K$ is a satellite knot. Set
 \[ L(K) = \max_{\mathcal{T}} \{L(\mathcal{T}), \; \mathcal{T} \text{ is a JSJ torus of } S^3_K\}, \]
 where the $L(\mathcal{T})$'s are given by Proposition \ref{prop:sxpiecesbound}. Let $|q|>L(K)$.
	
	Suppose there is an orientation-preserving homeomorphism $S^3_K(p/q) \cong S^3_{K'}(p/q)$ that does not carry the surgered piece $X$ of $S^3_{K}(p/q)$ to the surgered piece $X'$ of $S^3_{K'}(p/q)$. Then, as described in Section \ref{sec:scenario}, $X'$ is the image of a JSJ piece of $S^3_{K}(p/q)$ that is the outermost piece of the exterior of some knot $J$ such that $S^3_J \subset S^3_K$. Let $\mathcal{T} = \del S^3_{J}$. 
 
 Proposition \ref{prop:sxpiecesbound} implies that since $|q| > L(K) \geq L(\mathcal{T})$, the outermost piece of $S^3_J$ cannot be mapped to the surgered piece $X'$, a contradiction. Therefore, if $|q| > L(K)$, then any orientation-preserving homeomorphism $S^3_K(p/q) \cong S^3_{K'}(p/q)$ must send the surgered pieces one to another.
	\end{proof}

The proof of Proposition \ref{prop:sxpiecesbound} is divided into three cases: composing patterns, once or twice-iterated cables, and other patterns.

For the last two cases, we will need a simplified version of a theorem from Cooper and Lackenby, as well as some related lemmas.

\begin{thm}[{\cite[Theorem 4.1]{cooper}}]\label{thm:finitehyperb}
Let $M$ be a compact orientable 3-manifold with boundary a union of tori. Let $\epsilon > 0$. Then there are finitely many compact orientable hyperbolic 3-manifolds $X$ and slopes $\sigma$ on some component of $\del X$ such that $M \cong X(\sigma)$ and where the length of each slope $\sigma$ is at least $2\pi+\epsilon$, when measured using some horoball neighbourhood of the cusp of $X$ that is being filled.
\end{thm}

\begin{lemma}\label{lemma:hyplengthbound}
Let $Y$ be a hyperbolic JSJ piece of a knot exterior and let $\mathcal{L}_0$ be the cusp of $Y$ along which the trivial filling yields the exterior of an unlink. Let $l(p/q)$ be the length of the slope $p/q$ on $\mathcal{L}_0$, measured in a maximal horoball neighbourhood $N$ of $\mathcal{L}_0$. Then $l(p/q) \geq |q|/\sqrt{3}$.
\end{lemma}
\begin{proof}
By a geometric argument as in \cite[Lemma 2.1]{cooper} or \cite[Theorem 8.1]{agol_bounds_2000}, the lengths of two slopes $\sigma_1, \sigma_2$ on $\mathcal{L}_0$ satisfy $l(\sigma_1)l(\sigma_2) \geq Area(\del N) \cdot \Delta(\sigma_1, \sigma_2)$. 
By Theorem 1.2 of \cite{gabai2021hyperbolic}), $Area(\del N) \geq 2\sqrt{3}$.
By taking $\sigma_1 = p/q$ and $\sigma_2 = 1/0$, and by the 6-theorem (\cite{agol_bounds_2000}, \cite{lackenby_word_2000}), we get
\[l(p/q)\geq 2\sqrt{3} \cdot |q|/l(1/0)  \geq |q|/\sqrt{3}. \qedhere\]
\end{proof}
The next lemma follows the approach of \cite[Theorem 1.1, Case 2]{lack}.

\begin{lemma}\label{lemma:hypLY}
Let $Y$ be a JSJ piece of the exterior of a knot. There exists a constant $L(Y)$ with the following property. Let $Y'$ be a hyperbolic JSJ piece of the exterior of a knot, with boundary component $\mathcal{L}_0$ such that $Y'(\mathcal{L}_0; 1/0)$ is $S^3$ or the exterior of an unlink. If $Y'(\mathcal{L}_0; p/q) \cong Y$, then $|q|\leq L(Y)$.
\end{lemma}
\begin{proof}
Let $\epsilon = 1/15$. By Theorem \ref{thm:finitehyperb}, there are finitely many manifolds $\{X_j\}$ that are JSJ pieces of a knot exterior, and finitely many slopes $\{p_{j_i}/q_{j_i}\}$ of length at least $2\pi + 1/15$ (measured in a maximal horoball neighbourhood in $X_j$) such that $X_j(p_{j_i}/q_{j_i})$ is homeomorphic to $Y$. Set $L(Y) = \max \{|q_{j_i}|, 11\}$. If $|q|>L(Y)$, then by the previous lemma, $l(p/q) > 11/\sqrt{3} \geq 2\pi + 1/15$, but $p/q \notin \{p_{j_i}/q_{j_i}\}$. It follows that for any hyperbolic JSJ piece $Y'$ as in the statement, the filling $Y'(\mathcal{L}_0; p/q)$ cannot be homeomorphic to $Y$.
\end{proof}

\subsection{Composing pattern case}
We begin the proof of Proposition \ref{prop:sxpiecesbound} by considering the case of composing patterns.

\begin{lemma}\label{lemma:compositeknotext}
    Let $P$ be a composing pattern and $\mathcal{P} = \del \nu P \subset \del V_P$. If $|q|>1$, then the filling $V_P(\mathcal{P}; p/q)$ is not homeomorphic to a knot exterior.
\end{lemma}
\begin{proof}
    Let $Y \subset V_P$ be the composing space containing $\mathcal{P}$. Let $n+1$ be the number of boundary components of $Y$.
    
    If $n > 2$, then $Y(\mathcal{P};p/q)$ is Seifert fibred and has more than two boundary components (proof of Proposition \ref{prop:fillingY0}(2)), so it is not a JSJ piece of a knot exterior by Proposition \ref{prop:geomswitch}.

    Suppose now that $n=2$. Then $V_P = Y \cup S^3_{K_1}$ for some knot $K_1$.
The filling $Y(\mathcal{P};p/q)$ is Seifert fibred (proof of Proposition \ref{prop:fillingY0}(2)), and on the JSJ torus $\mathcal{T}_1 = \del S^3_{K_1}$ of $V_P(\mathcal{P};p/q)$, a regular fibre of $Y(\mathcal{P};p/q)$ has meridional slope. 

By Lemma \ref{prop:geomswitch}, if $J'$ is a knot such that $S^3_{J'}$ has the same JSJ pieces in its decomposition as $V_P(\mathcal{P};p/q)$, then $J'$ must be a cable of $K_1$. By the knot complement theorem, if $V_P(\mathcal{P};p/q)$ were homeomorphic to $S^3_{J'}$, then the meridian on $\mathcal{T}_1$ would be mapped to the meridian on $\mathcal{T}_1' = \del S^3_{K_1} \subset S^3_{J'}$. Further, regular fibres of $Y(\mathcal{P};p/q)$ would be mapped to regular fibres of the outermost cable space of $S^3_{J'}$. However, a regular fibre of the outermost cable space of $S^3_{J'}$ does not have meridional slope on $\mathcal{T}_1'$, and a cable space possesses a unique Seifert fibred structure. Hence, $V_P(\mathcal{P};p/q)$ cannot be homeomorphic to the exterior of a knot.
\end{proof}

\begin{prop}\label{prop:compositeswitch}
    Let $K = K_1 \# K_2 \# \ldots \# K_n$ be a composite knot, where the $K_i$'s are prime for each $i=1, \ldots, n$. Let $\mathcal{T}$ be a JSJ torus of $S^3_K$ that decomposes $K$ into $P$, a composing pattern, and $K_i$ for some $i \in \{1, \ldots, n\}$. Suppose there is an orientation-preserving homeomorphism $S^3_K(p/q) \cong S^3_{K'}(p/q)$ for some knot $K'$ where $|q|>2$. Then the homeomorphism does not map the outermost piece of $S^3_{K_i}$ to the surgered piece of $S^3_{K'}(p/q)$.
\end{prop}
\begin{proof}

Let $\mathcal{P} = \del \nu P \subset \del V_P$. If the homeomorphism maps the outermost piece of $S^3_{K_i}$ to the surgered piece of $S^3_{K'}(p/q)$, then $V_P(\mathcal{P}; p/q)$ is homeomorphic to a knot exterior by the discussion of Section \ref{sec:scenario} and Figure \ref{fig:switchinghomeo}. By Lemma \ref{lemma:compositeknotext}, this cannot happen if $|q|>2$.
\end{proof}

\begin{proof}[Proof of Proposition \ref{prop:sxpiecesbound}]\renewcommand{\qedsymbol}{}
    Let $\mathcal{T}$ be a JSJ torus of $S^3_K$ that decomposes $K$ into $P$ and $J$. Suppose there exists a knot $K'$ such that there is an orientation-preserving homeomorphism $S^3_K(p/q) \cong S^3_{K'}(p/q)$. By Proposition \ref{prop:compositeswitch}, if $P$ is a composing pattern, we may take $L(\mathcal{T}) = 2$. 
\end{proof}
\subsection{Cable and twice-iterated cable case}\label{sec:surgeredcable}

We proceed with the case when $P$ is a once or twice-iterated cable.

\begin{proof}[Proof of Proposition \ref{prop:sxpiecesbound} (continued)]\renewcommand{\qedsymbol}{}
We now suppose that $P$ is a cable $C_{r_1,s_1}$ or a twice-iterated cable $C_{r_2,s_2}(C_{r_1,s_1})$. Recall that the JSJ torus $\mathcal{T}$ decomposes $K$ into $P$ and a knot $J$. Let $Y$ be the outermost piece of $S^3_J$. Let $Y'$ be the JSJ piece of $S^3_{K'}$ such that the surgered piece of $S^3_{K'}(p/q)$ is $X' = Y'(p/(q(t')^2), t' \geq 1$ (Proposition \ref{prop:geomsxpiece}). 

Suppose the homeomorphism $S^3_K(p/q) \cong S^3_{K'}(p/q)$ carries the outermost piece $Y$ of $S^3_J$ to the surgered piece of $S^3_{K'}(p/q)$, as described in Section \ref{sec:scenario}. We look at each possibility for $Y$ given by Proposition \ref{prop:geomswitch}.


If $Y$ is hyperbolic, then $Y'$ is also hyperbolic if $|q|>2$, according to Proposition \ref{prop:geomswitch}(2). By Lemma \ref{lemma:hypLY}, there is a constant $L(Y)$ such that $|q| \leq |q(t')^2| \leq L(Y)$.


If $Y$ is an $(r,s)$-cable space, then $Y'$ is a composing space by Proposition \ref{prop:geomswitch}(3) and $K'$ is a composite knot. Using the notation in Figure \ref{fig:switchinghomeo}, $\mathcal{T}'$ separates $K'$ into a composing pattern $P'$ and some companion knot $J'$. By the discussion of Section \ref{sec:scenario}, $V_{P'}(\mathcal{P}'; p/q)$ is homeomorphic to $S^3_J$, but this contradicts Lemma \ref{lemma:compositeknotext} applied to $P'$ when $|q|>2$.

If $Y$ is the exterior of a torus knot $T_{a,b}, |a|>|b|>1$, then $Y'$ is an $(r',s')$-cable space by Proposition \ref{prop:geomswitch}(4). Since the orders of exceptional fibres in $S^3_{T_{a,b}}$ and $X'$ coincide, we have without loss of generality
	\[|a| = |q(t')^2r's'-p|. \tag{1}\]

Recall that the JSJ torus $\mathcal{T} \subset \del Y$ is mapped by the homeomorphism $S^3_K(p/q) \cong S^3_{K'}(p/q)$ to the JSJ torus $\mathcal{T}' \subset \del X'$, in the notation of Figure \ref{fig:switchinghomeo}. As distances are preserved between slopes that are carried one to another, we have the following equality by comparing Table \ref{table:distslopesCT} (regular fibre slopes) and the last row of Table \ref{table:distslopes} (longitudinal slopes) from Section \ref{sec:distslopes}:
\[|r_1| = \Delta\left( \frac{r_1}{s_1}, \frac{0}{1}\right) = \Delta\left( \frac{k}{1}, \frac{p}{q(t's')^2}\right)=|q(t's')^2k-p|.\]
Combining this with equation (1) yields
\[|q(t')^2| \cdot |(s')^2k-r's'| = | r_1 \pm a|.\]
Since $r',s'\neq 0$ are coprime, we have $(s')^2k-r's' \neq 0$. This implies that $|q| \leq |r_1|+|a|$.

Summing up, suppose $\mathcal{T}$ decomposes $K$ into $P$ and $J$, where $P = C_{r_1,s_1}$ or $C_{r_2,s_2}(C_{r_1,s_1})$. Denoting the outermost JSJ piece of $S^3_J$ by $Y$, we let
 \[ L(\mathcal{T}) = 
 \begin{cases}
 L(Y) \text{ from Lemma \ref{lemma:hypLY}}& \text{if } Y \text{ is hyperbolic},  \\
 2 & \text{if } Y \text{ is an } (r,s)\text{-cable space}, \\
 |r_1|+|a| & \text{if } Y \text{ is the exterior of a torus knot } T_{a,b}.
 \end{cases}\]
 Then if $|q| > L(\mathcal{T})$, and if there exists and orientation-preserving homeomorphism $S^3_K(p/q) \cong S^3_{K'}(p/q)$ for some knot $K'$, the outermost piece $Y$ of $S^3_J$ is not carried to the surgered piece of $S^3_{K'}(p/q)$.
 \end{proof}

\subsection{Other pattern case}
To conclude the proof of Proposition \ref{prop:sxpiecesbound}, it remains to study patterns that are neither composing patterns nor once or twice-iterated cables. We will be using the Cyclic surgery theorem by Culler, Gordon, Luecke and Shalen.

\begin{thm}[{\cite[Cyclic surgery theorem]{culler_dehn_1987}}]\label{thm:cyclic}
Let $M$ be a compact, connected, irreducible, orientable 3-manifold such
that $\del M$ is a torus. Suppose that $M$ is not a Seifert fibre space. If
$\pi_1 (M(\sigma_1))$ and $\pi_1(M(\sigma_2))$ are cyclic, then $\Delta(\sigma_1, \sigma_2) \leq 1$.
\end{thm}

We want to apply this theorem to the case where $M$ is a filling of a pattern space. To do so, we must show that this filling is not a Seifert fibre space. We will need the following homological lemma about fillings of composing spaces.

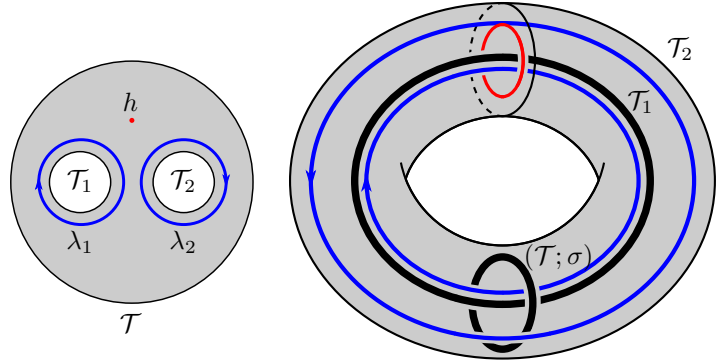
\begin{figure}
    \centering
    \begin{tikzpicture}
    \node at (0,0) {\includegraphics[scale=0.37,trim = 0pt -57pt 0pt 0pt]{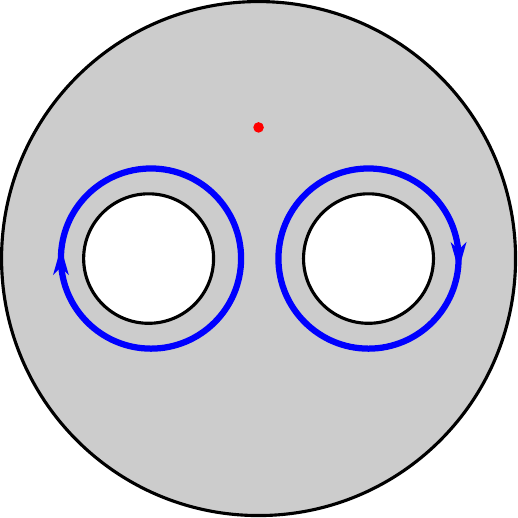} \quad \includegraphics[scale=0.53]{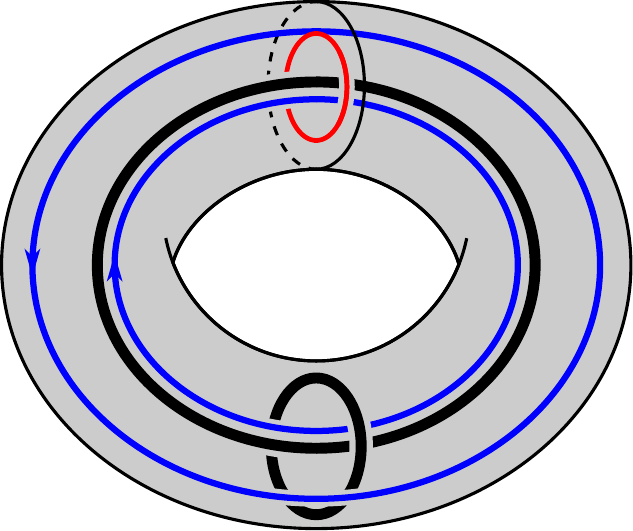}};

    \node at (-107pt,0) {$\mathcal{T}_1$};
    \node at (-68pt,0) {$\mathcal{T}_2$};
    \node at (-88pt,30pt) {$h$};
    \node at (-107pt,-23pt) {$ \lambda_1$};
    \node at (-68pt,-23pt) {$\lambda_2$};    
    \node at (-88pt,-54pt) {$\mathcal{T}$};
    
    \node at (74pt,-29pt) {$(\mathcal{T}; \sigma)$};
    \node at (105pt,30pt) {$\mathcal{T}_1$};
    \node at (120pt,50pt) {$\mathcal{T}_2$};
    
    \end{tikzpicture}
    
    \caption{Left: Generators of $H_1(Y; \Z)$ on the base orbifold of $Y$; Right: Generators of $H_1(Y(\mathcal{T}; \sigma); \Z)$ where $Y$ is depicted as a link exterior in a solid torus bounded by $\mathcal{T}_2$}
    \label{fig:genH1composing}
\end{figure}

\begin{lemma}\label{lemma:H1composing} Let $Y$ be a composing space with three boundary components $\mathcal{T}, \mathcal{T}_1, \mathcal{T}_2$. Denote by $h$ the slope of a regular fibre on each boundary component of $Y$. Suppose $\sigma_1$ and $\sigma_2$ are slopes on $\mathcal{T}_1$ and $\mathcal{T}_2$ respectively, that are homologous in $Y(\mathcal{T}; \sigma)$ for some surgery slope $\sigma$ on $\mathcal{T}$. Then $\Delta(h, \sigma_1) = \Delta(h, \sigma_2) = k\Delta(h, \sigma)$ for some $k \in \Z$.
\end{lemma}
\begin{proof}
    There are slopes $\lambda_1$ and $\lambda_2$ on $\mathcal{T}_1$ and $\mathcal{T}_2$ respectively such that $\{h, \lambda_i\}$ generates $H_1(\mathcal{T}_i; \Z)$, $i=1,2$, and $\{h, \lambda_2-\lambda_1\}$ generates $H_1(\mathcal{T}; \Z)$ (Figure \ref{fig:genH1composing}). Further, the images induced by inclusion of $h, \lambda_1, \lambda_2$ into $Y$ generate $H_1(Y; \Z)$. Write 
    \begin{align*}
      \sigma &= mh + n(\lambda_2-\lambda_1),\\
      \sigma_1 &= a_1h + b_1\lambda_1, \\
      \sigma_2 &= a_2h + b_2\lambda_2.
    \end{align*}
    The $\sigma$-surgery along $\mathcal{T}$ adds the relation $mh + n(\lambda_2-\lambda_1)$ in $H_1(Y(\mathcal{T}; \sigma); \Z)$. Hence, if $\sigma_1$ and $\sigma_2$ are homologous in $Y(\mathcal{T}; \sigma)$, then
    \[(a_1h + b_1\lambda_1) + k(mh + n(\lambda_2-\lambda_1)) = a_2h + b_2\lambda_2,\]
    for some $k \in \Z$. This implies that $b_1=b_2=kn$, giving us
    \[\Delta(h, \sigma_1) = \Delta(h, \sigma_2) = kn = k\Delta(h, \sigma).\qedhere\]
\end{proof}

\begin{prop}\label{prop:notSeifert fibre space}
Suppose $P$ is a pattern that is neither a composing pattern, nor a cable $C_{r_1,s_1}$ or a twice-iterated cable $C_{r_2,s_2}(C_{r_1,s_1})$. Let $\mathcal{T}$ be the boundary component of $V_P$ that is not $\del \nu P$. Then the filling $V_P(\mathcal{T};m/n)$ is not a Seifert fibre space for $|m|$ sufficiently large.
\end{prop}
\begin{proof}
The pattern space $V_P$ admits a JSJ decomposition $V_1 \cup V_2 \cup \ldots \cup V_k$ where $V_1$ is the JSJ piece that contains $\mathcal{T}$. Recall from the discussion following Definition \ref{def:Tdecomposes} that we express slopes along $\mathcal{T}$ in the coordinates given by slopes $\mu$ and $\lambda$ such that gluing the exterior of a knot $J$ by respectively identifying $\mu$ and $\lambda$ to the meridian and longitude of $J$ yields the exterior of the knot $P(J)$.

If $V_1$ is hyperbolic, there are only finitely many slopes $m/n$ such that $V_1(\mathcal{T}; m/n)$ is not hyperbolic. So $V_1(\mathcal{T};m/n)$ is hyperbolic for $|m|$ sufficiently large, and $V_P(\mathcal{T}; m/n)$ is not a Seifert fibre space.

If $V_1$ is a Seifert fibre space, then by Theorem \ref{thm:budney}, $V_1$ is either a composing space or an $(r_1,s_1)$-cable space. 

Suppose $V_1$ is a composing space. For $V_P(\mathcal{T}; m/n)$ to be Seifert fibered, the JSJ pieces adjacent to $V_1$ in $V_P$ must be Seifert fibred and $V_1(\mathcal{T};m/n)$ must admit a Seifert fibred structure that differs from the one inherited by the fibration on $V_1$. The only such possibility is if $V_1(\mathcal{T};m/n)$ is a trivial $I$-bundle over the torus. Recall that the regular fibres of $V_1$ have meridional slopes on each boundary component of $V_1$. Let $\mathcal{T}_1$ and $\mathcal{T}_2$ be the boundary components of $V_1(\mathcal{T};m/n)$. As regular fibres are homologous in $V_1(\mathcal{T};m/n)$, Lemma \ref{lemma:H1composing} says the distance on $\mathcal{T}_1$ between the meridian of $\mathcal{T}_1$ and a regular fibre of $V_1(\mathcal{T};m/n)$ is equal to the distance on $\mathcal{T}_2$ between the meridian of $\mathcal{T}_2$ and a regular fibre of $V_1(\mathcal{T};m/n)$.

Suppose the pieces adjacent to $V_1$ in $V_P$ are Seifert fibred. Note that since $P$ is not a composing pattern, $V_1(\mathcal{T};m/n)$ shares a boundary component, say $\mathcal{T}_1$, with a cable space $V_2$ whose regular fibre has non-integral slope on $\mathcal{T}_1$. The other boundary component $\mathcal{T}_2$ of $V_1(\mathcal{T};m/n)$ is shared with a torus knot exterior or a cable space $V_3$, whose regular fibre has integral slope on $\mathcal{T}_2$. By Lemma \ref{lemma:H1composing}, the Seifert fibred structure of $V_1(\mathcal{T};m/n)$ cannot extend across both $V_2$ and $V_3$, so $V_P(\mathcal{T};m/n)$ is not Seifert fibred.

Suppose now that $V_1$ is an $(r_1, s_1)$-cable space. Let $V_2$ be the JSJ piece of $V_P$ that shares a boundary component $\mathcal{T}_1$ with $V_1$. 

Suppose that $\mathcal{T}_1$ remains incompressible in $V_P(\mathcal{T}; m/n)$. The pattern space $V_P$ is either the union of $V_1$ with a hyperbolic $V_2$, or it decomposes into at least three JSJ pieces. In the first case, $V_P(\mathcal{T}; m/n)$ is clearly not Seifert fibred. In the second case, a Seifert fibred structure on $V_1(\mathcal{T}; m/n)$ might extend across a Seifert fibred structure on $V_2$. However, a JSJ piece of a knot exterior admits a unique Seifert fibred structure, so the structure on $V_2$ does not extend across the other JSJ pieces of $V_P(\mathcal{T}; m/n)$.

Suppose now that the torus $\mathcal{T}_1$ is compressed in $V_1(\mathcal{T}; m/n)$. On $\mathcal{T}_1$ and $\mathcal{T}$, the regular fibres of $V_1$ have respective slopes $r_1s_1/1$ and $r_1/s_1$. By a similar reasoning as that of Proposition \ref{prop:fillingY0}, cases (3) and (4), we have $|ms_1-r_1n|=1$. For homological reasons (analogous to Lemma \ref{lemma:nullhom}), the filling $V_P(\mathcal{T};m/n)$ is homeomorphic to $(V_P \setminus V_1)(\mathcal{T}_1; ms_1^2/n)$. 

If $V_2$ is hyperbolic or a composing space, we iterate the argument previously given for $V_1$.

If $V_2$ is an $(r_2, s_2)$-cable space, let $\mathcal{T}_2$ be its boundary component that is not $\mathcal{T}_1$. Let $V_3$ be the JSJ piece of $V_P$ such that $V_3 \cap V_2 = \mathcal{T}_2$. The Seifert fibred structure on $V_2(\mathcal{T}_1; ms_1^2/n)$ might extend across $V_3$ only if $V_2(\mathcal{T}_1; ms_1^2/n)$ is a solid torus or a twisted $I$-bundle over the Klein bottle. This occurs when $|ms_1^2s_2-r_2n| = 1$ or $2$. Combining this with the fact that $|ms_1-r_1n|=1$, we have that $(m,n)$ must be a solution to the system
\[\left(\begin{matrix} s_1 & -r_1 \\ s_1^2s_2 & -r_2 \end{matrix}\right) \left(\begin{matrix} m  \\ n \end{matrix}\right) 
= \left(\begin{matrix} \pm 1  \\ \pm 1 \end{matrix}\right) \text{ or } \left(\begin{matrix} \pm 1  \\ \pm 2 \end{matrix}\right).\]
As $r_2$ and $s_2$ are coprime, \[\det \left(\begin{matrix} s_1 & -r_1 \\ s_1^2s_2 & -r_2 \end{matrix}\right) \neq 0.\]
Therefore, there are only finitely many slopes $m/n$ such that the Seifert fibred structure on $V_2(\mathcal{T}_1; ms_1^2/n)$ extends across $V_3$.

Consequently, $V_P(\mathcal{T};m/n)$ is not a Seifert fibre space for $|m|$ sufficiently large.
\end{proof}

\begin{proof}[Proof of Proposition \ref{prop:sxpiecesbound} (continued)]
Recall that $\mathcal{T}$ is a JSJ torus of $S^3_K$ that decomposes $K$ into $P$ and $J$. Suppose there exists a knot $K'$ such that there is an orientation-preserving homeomorphism $S^3_K(p/q) \cong S^3_{K'}(p/q)$ where $|q|>2$.

We now suppose that $P$ is neither a composing pattern nor a once or twice-iterated cable. Suppose the homeomorphism $S^3_K(p/q) \cong S^3_{K'}(p/q)$ carries the outermost piece $Y$ of $S^3_J$ to the surgered piece of $S^3_{K'}(p/q)$, as described in the scenario of Section \ref{sec:scenario}. By Lemma \ref{lemma:meridian}, there is a slope $x/y = q(w')^2/y$ on $\mathcal{T}$ that is the meridian of $V_P(\mathcal{P}; p/q)$ seen as a knot exterior, where $w'$ is the winding number of $P'$. See Figure \ref{fig:switchinghomeo}.

By Lemma \ref{prop:notSeifert fibre space}, there exists a bound $L(P)$ such that if $|m|>L(P)$, then $V_P(\mathcal{T}; m/n)$ is not a Seifert fibre space. If $w'\neq 0$, suppose that $|q|>L(P)$. The inequality $|x|=|q(w')^2| > |q| > L(P)$ implies that $V_P(\mathcal{T}; x/y)$ is not a Seifert fibre space. On one hand, the filling of $V_P(\mathcal{P}; p/q)$ along the meridian $x/y$ is the trivial filling $V_P(\mathcal{P}, \mathcal{T}; p/q, x/y) \cong S^3$. On the other hand, $V_P(\mathcal{P}; 1/0)$ is the trivial filling of the pattern $P$, so it is homeomorphic to a solid torus. Consequently, $V_P(\mathcal{P}, \mathcal{T}; 1/0, x/y)$ is homeomorphic to the lens space $L_{y,x}$. We obtain that both the $p/q$ and $1/0$-fillings of the non-Seifert fibred manifold $V_P(\mathcal{T}; x/y)$ yield manifolds with cyclic fundamental groups. By the Cyclic surgery theorem (Theorem \ref{thm:cyclic}), we have $|q| = \Delta(p/q, 1/0) \leq 1$, which contradicts $|q|>2$.

Suppose now that $w'=0$. If the surgered piece $X'$ of $S^3_{K'}(p/q)$ were Seifert fibred, then it would be a filling of either a cable space or a composing space (Proposition \ref{prop:geomswitch}). In both cases, $w'$ would be non-zero, a contradiction. Therefore, $X'$ is hyperbolic. As $|q|>2$, $X'$ is the $p/(q(t')^2)$-filling of a hyperbolic JSJ piece of $S^3_{K'}$, $t' \geq 1$ (Proposition \ref{prop:geomsxpiece}(1)). By Lemma \ref{lemma:hypLY}, there exists a constant $L(Y)$ such that $X'$ is homeomorphic to the outermost piece $Y$ of $S^3_J$ only if $|q|\leq L(Y)$.

Setting $L(\mathcal{T})= \max \{ L(P), L(Y) \}$ gives the desired bound.

This completes the proof of Proposition \ref{prop:sxpiecesbound}. 
\end{proof}

  \section{Proof of Theorem 1}\label{sec:proof}

Proposition \ref{prop:sxpiece} tells us that if $|q|$ is sufficiently large, an orientation-preserving homeomorphism $S^3_K(p/q) \cong S^3_{K'}(p/q)$ restricts to a homeomorphism between the surgered pieces of $S^3_K(p/q)$ and $S^3_{K'}(p/q)$. By the knot complement theorem, this homeomorphism preserves the slopes on the boundary of the surgered pieces. To complete the proof of Theorem \ref{thm:main}, we must show that it further restricts to the JSJ pieces of $S^3_K$ and $S^3_{K'}$ that were filled to produce the surgered pieces.

First, we need the following intermediate results.

\begin{prop}\label{prop:std}
Let $K$ and $K'$ be knots such that there exists an orientation-preserving homeomorphism $S^3_K(p/q) \cong S^3_{K'}(p'/q')$. If the core of the surgery solid torus in $S^3_K(p/q)$ is mapped to the core of the surgery solid torus in $S^3_{K'}(p'/q')$, then $K = K'$.
\end{prop}
\begin{proof}
Let $v$ and $v'$ be the cores of the surgery solid tori of $S^3_K(p/q)$ and $S^3_{K'}(p'/q')$ respectively. Since $v$ is sent to $v'$ by the homeomorphism, the neighbourhoods $\nu(v)$ and $\nu(v')$ are also sent one to another by the homeomorphism. Therefore, $S^3_K(p/q) \setminus \operatorname{int}(\nu(v)) \cong S^3 \setminus \operatorname{int}(\nu K)$ is homeomorphic to $S^3_{K'}(p'/q') \setminus \operatorname{int}(\nu(v')) \cong S^3 \setminus \operatorname{int}(\nu K')$, which implies that $K=K'$ by the knot complement theorem.
\end{proof}

\begin{lemma}\label{lemma:thelemma}
If $q$, $p$, $r$, $s$, $r'$, $s'$ are integers such that $|q|>2$ and $|qrs-p|=|qr's'-p|=1$, then $rs = r's'$.
\end{lemma}
\begin{proof}
We have $|q(rs-r's')| = 0$ or $2$. But $|q|>2$, so $|q(rs-r's')| = 0$ and $rs=r's'$.
\end{proof}

Let $K$ be a non-trivial knot, and suppose there is a knot $K'$ such that there exists an orientation-preserving homeomorphism $S^3_K(p/q) \cong S^3_{K'}(p/q)$. Let $X$, $X'$ be the surgered pieces of $S^3_{K}(p/q)$, $S^3_{K'}(p/q)$ respectively. Let $Y$, $Y'$ be the JSJ pieces of $S^3_{K}$, $S^3_{K'}$ such that $X=Y(p/(qt^2))$ and $X'=Y'(p/(q(t')^2))$, for some $t,t' \geq 1$ (Proposition \ref{prop:geomsxpiece}). 

We now assume, by Proposition \ref{prop:sxpiece}, that $|q|$ is large enough such that $X$ and $X'$ are sent one to another by the homeomorphism $S^3_K(p/q) \cong S^3_{K'}(p/q)$. We will study each possibility listed in Theorem \ref{thm:budney} for $Y$ and show in each case that $K=K'$ for $|q|$ sufficiently large.

\subsection{Exterior of a torus knot}\label{sec:extTorusknot}

Suppose $Y$ is the exterior of a torus knot. In this case, $K$ is a torus knot or a cable of a torus knot (Proposition \ref{prop:geomsxpiece}). By McCoy, if $K$ is a torus knot, we have that $K = K'$ for $|q|$ sufficiently large (\cite{mcc}). So suppose $K$ is a cable knot $C_{r,s}(T_{a,b})$ such that $|qrs-p|=1$.

By Corollary \ref{cor:sxtype}, $Y'$ is also the exterior of a torus knot if $|q|>8$. Therefore, $K'$ is either a torus knot or a cable of a torus knot by Proposition \ref{prop:JSJsurgery}.

We have the following corollary of a proposition from McCoy.
\begin{prop}[{\cite[Proposition 1.5]{mcc}}]
If an $(r,s)$-cable of a torus knot shares a $p/q$-surgery with a torus knot where $|q|>1$, then $|q|=s$.
\end{prop}
It follows that the cable $K=C_{r,s}(T_{a,b})$ cannot share a $p/q$ surgery with a torus knot when $|q| > s$. Hence, if $|q| > s$ and $8$, $K'$ is a cable of a torus knot $C_{r', s'}(T_{c,d})$ where $|qr's'-p|=1$. We thus have an homeomorphism
\[S^3_{T_{a,b}}(p/(qs^2)) \cong S^3_{T_{c,d}}(p/(q(s')^2))\]
by Proposition \ref{prop:cablepqs2}. This gives the homeomorphism of base orbifolds
\[S^2(|a|, |b|, |qs^2ab-p|) \cong S^2(|c|, |d|, |q(s')^2cd-p|).\]

Comparing orders of cone points, without loss of generality, assume that $|b| = |d|$ and $|a| = |q(s')^2cd-p|$. By combining this with $|qrs-p|=1$, we find
\[|q| \cdot |(s')^2cd-rs)| = | a \pm 1|.\]
The right-hand side is a non-zero integer since $|a|>1$, which implies that $|q| \leq |a|+1$. 

Consequently, if $|q| > |a|+1$, the homeomorphism $S^3_{T_{a,b}}(p/(qs^2)) \cong S^3_{T_{c,d}}(p/(q(s')^2))$ sends the core of the surgery solid torus in $S^3_{T_{a,b}}(p/(qs^2))$ of order $|qs^2ab-p|$ to the core of the surgery solid torus in $S^3_{T_{c,d}}(p/(q(s')^2))$ of order $|q(s')^2cd-p|$. By Proposition \ref{prop:std}, we obtain that $T_{a,b} = T_{c,d}$. Furthermore, the equality of orders yields $qs^2ab-p=\pm (q(s')^2cd-p)$, but since $|q|>1$ and $p$ and $q$ are coprime, the only possibility is $qs^2ab-p=q(s')^2cd-p$, which in turn gives $s=s'$.  By Lemma \ref{lemma:thelemma}, since $|q|>|a|+1 > 2$, we have $C_{r,s} = C_{r', s'}$. Hence, $C_{r,s}(T_{a,b}) = C_{r', s'}(T_{c,d})$, that is, $K=K'$, as desired.

\subsection{Composing space}\label{sec:theorcompsing}

Suppose $Y$ is a composing space. By Corollary \ref{cor:sxtype}, $Y'$ is also a composing space if $|q|>2$. By Proposition \ref{prop:geomsxpiece}(2), $X$ and $X'$ are Seifert fibred, each with one exceptional fibre of order $|qt^2|$ and $|q(t')^2|$ respectively. These exceptional fibres correspond to the cores of the surgery solid tori in $X$ and $X'$. 

Since $X \cong X'$, the unique exceptional fibre of $X$ is sent to the unique exceptional fibre of $X'$ by the homeomorphism $S^3_{K}(p/q) \cong S^3_{K'}(p/q)$. This implies that $t = t'$. If $t = 1$, then these exceptional fibres are precisely the cores of the surgery solid tori of $S^3_{K}(p/q)$ and $S^3_{K'}(p/q)$. Then $K' = K$, by Proposition \ref{prop:std}. If $t > 1$, by Proposition \ref{prop:JSJsurgery}, $t$ is the winding number of cable patterns $C_{r,t}$ and $C_{r',t}$ such that $K=C_{r,t}(J)$ and $K'=C_{r',t}(J')$, where $J$ and $J'$ are composite knots. By Proposition \ref{prop:cablepqs2}, we have $S^3_{J}(p/(qt^2)) \cong S^3_{J'}(p/(qt^2))$ and $J = J'$ by Proposition \ref{prop:std}. Since $|qrt-p|=1 $ and $|q r' t-p|=1$, Lemma \ref{lemma:thelemma} tells us that $r=r'$, and we conclude that $K = K'$.

\subsection{Exterior of a hyperbolic link}
Suppose $Y$ is the exterior of a hyperbolic knot or link. By Corollary \ref{cor:sxtype}, $Y'$ is also the exterior of a hyperbolic knot or link if $|q|>8$. Recall that $Y(p/(qt^2)) \cong Y'(p/(q(t')^2))$. We will apply the following theorem by Lackenby and use the arguments in his proof of \cite[Case 2 of Theorem 1.1]{lack}.

\begin{thm}[{\cite[Theorem 3.1]{lack}}]\label{thm:CKhyperb}
	Let $M$ be $S^3$ or the exterior of the unknot or unlink in $S^3$, and let $K$ be a hyperbolic knot in $M$. Let $M_K = M \setminus \operatorname{int}(\nu K)$. There exists a constant $C(K)$ with the following property. If $M_K(\sigma) \cong M_{K'}(\sigma')$ for some hyperbolic knot $K'$ in $M$ and some $\sigma'$ such that $\Delta(\sigma',\mu') > C(K)$, where $\mu'$ is the slope that bounds a disc in $\nu K'$, and if the homeomorphism restricted to the boundary of $M$ is the identity, then $(M, K) \cong (M, K')$ and $\sigma = \sigma'$.
\end{thm}

\begin{lemma}\label{lemma:applyLack}
    For $|q|$ sufficiently large, $Y \cong Y'$ and $t=t'$.
\end{lemma}
\begin{proof}
Let $n+1$ be the number of boundary components of $Y$ and $Y'$. If $n=0$, let $M$ be $S^3$. If $n\geq 1$, let $M$ be the exterior of the unlink with $n$ components. By Theorem \ref{thm:budney}, $Y$ and $Y'$ are respectively homeomorphic to exteriors of hyperbolic knots $H$ and $H'$ in $M$.

By the argument in the last paragraph of \cite[proof of Theorem 1.1, p.13]{lack}, there is a knot $H''$ in $M$ such that $(M,H') \cong (M,H'')$ and there exists an homeomorphism $M_{H}(p/(qt^2)) \cong M_{H''}(p/(q(t')^2))$ which is the identity when restricted to the boundary of $M$. Let $C(H)$ be the constant given by Theorem \ref{thm:CKhyperb} for $H$. If $|q| > C(H)$, then $|q(t')^2|>C(H)$. By Theorem \ref{thm:CKhyperb}, $(M,H) \cong (M,H'')$ and $p/(qt^2)$ = $p/(q(t')^2)$. Therefore, $(M, H) \cong (M, H')$, so $Y \cong Y'$, and $t=t'$.
\end{proof}

In $S^3_K$ and $S^3_{K'}$ respectively, the JSJ pieces $Y, Y'$ are the outermost pieces of the exteriors of knots $J, J'$. We thus have $S^3_J(p/(qt^2)) \cong S^3_{J'}(p/(qt^2))$ which restricts to $Y(p/(qt^2)) \cong Y'(p/(qt^2))$. This is precisely the scenario of \cite[Case 2 of Theorem 1.1]{lack}. We adapt the relevant parts of its proof using our notation to conclude the case when $Y$ is hyperbolic.

\begin{prop}
    Let $C(H)$ be as in the proof of Lemma \ref{lemma:applyLack}. If $|q|> \max\{8, C(H)\}$, then $K=K'$.
\end{prop}
\begin{proof}
The homeomorphism $(M, H) \cong (M, H')$ from the Lemma \ref{lemma:applyLack} gives a homeomorphism $h : S^3 \setminus (S^3_J \setminus \operatorname{int}(Y)) \rightarrow S^3 \setminus (S^3_{J'} \setminus \operatorname{int}(Y'))$ that sends $J$ to $J'$.

Further, $S^3_K(p/q) \cong S^3_{K'}(p/q)$ restricts to a homeomorphism $S^3_J \setminus \operatorname{int}(Y) \cong S^3_{J'} \setminus \operatorname{int}(Y')$ which agrees with $h$ on the boundary. We can thus extend it to a homeomorphism $(S^3, J) \cong  (S^3, J')$. Hence, $J = J'$. If $t = 1$, then $K=J$ and $K'=J'$ and we are done. If $t>1$, then $K$ and $K'$ are cables of $J = J'$. By Lemma \ref{lemma:thelemma}, $K = K'$.
\end{proof}

\subsection{Cable space}\label{sec:theorcable}
Suppose $Y$ is an $(r_1, s_1)$-cable space. By Corollary \ref{cor:sxtype}, $Y'$ is also a cable space if $|q|>2$. Therefore, $K$ and $K'$ are once or twice-iterated cables of knots $J$ and $J'$ respectively. Let us write
\[K = \begin{cases}
C_{r_1,s_1}(J) &\text{if }t=1 \\
C_{r_2,s_2}(C_{r_1,s_1}(J)) &\text{if }t>1
\end{cases} 
, \quad
K' = \begin{cases}
C_{r_1',s_1'}(J') &\text{if }t'=1 \\
C_{r_2',s_2'}(C_{r_1',s_1'}(J')) &\text{if }t'>1
\end{cases}.\]

\begin{prop} If $|q|>2$, then $K = K'$. That is:
	\begin{enumerate}[(i)]  
\setlength{\itemsep}{5pt}
  \setlength{\parskip}{0pt}
  \setlength{\parsep}{0pt}
		\item $J = J'$;		
		\item $C_{r_1,s_1} = C_{r_1',s_1'}$;
		\item $t=t'$, and $C_{r_2,s_2} = C_{r_2',s_2'}$ if $t>1$.
	\end{enumerate}
\end{prop}

\begin{proof}
	Since $S^3_{K}(p/q) \setminus X \cong S^3_J$ and $S^3_{K'}(p/q) \setminus X' \cong S^3_{J'}$, the homeomorphism $S^3_{K}(p/q) \cong S^3_{K'}(p/q)$ restricts to a homeomorphism between $S^3_{J}$ and $S^3_{J'}$. This implies (i) by the knot complement theorem.
	
	By (i), the meridians and longitudes forming the bases of $H_1(\del S^3_J;\Z)$ and $H_1(\del S^3_{J'};\Z)$ are respectively sent one to another.
The regular fibres of $X$ and $X'$ have respective slopes $r_1/s_1$ and $r_1'/s_1'$ on $\del X$ and $\del X'$, and both Seifert fibred structures have base orbifold a disc with two cone points. If a given oriented manifold admits a Seifert fibration with base orbifold a disc and two cone points, then there is no other Seifert fibration on this manifold with the same orbifold structure. It follows that the slopes $r_1/s_1$ and $r_1'/s_1'$ are equal. Hence, $C_{r_1,s_1} = C_{r_1',s_1'}$ showing (ii).

The longitudes of $X$ and $X'$ coincide and have respective slopes $p/(q(ts_1)^2)$ and $p/(q(t's_1')^2)$, so $q(ts_1)^2 = q(t's_1')^2$. Since $s_1 = s_1'$ by (ii), we get the equality $t = t'$. If $t>1$, then $t = s_2 = s_2'$, and $C_{r_2,s_2} = C_{r_2',s_2'}$ by Lemma \ref{lemma:thelemma}, which proves (iii).
\end{proof}

This concludes the proof of Theorem \ref{thm:main}. 
  \section{Characterizing slopes for cables with only Seifert fibred pieces}\label{sec:theorexplicitSFS}

For some specific families of satellite knots, an explicit bound for $|q|$ that realizes Theorem \ref{thm:main} can be expressed. The following result is obtained from the treatment of Seifert fibred JSJ pieces throughout Sections \ref{sec:sxpiece} and \ref{sec:proof}.

\begin{thm}\label{thm:SFSbounds}
Let $K$ be a cable knot with an exterior consisting solely of Seifert fibred JSJ pieces. A slope $p/q$ is characterizing for $K$ if:
    \begin{enumerate}[(i)]
        \item $|q|>2$ and $K$ is not an $n$-times iterated cable of a torus knot, $n \geq 1$;
        \item $|q|>|r_1| + |a|$ and $K$ is an $n$-times iterated cable of $C_{r_1, s_1}(T_{a, b}), |a|>|b|>1, n \geq 1$;
        \item $|q|> \max\{8, s_1, |r_1| + |a|\}$ and $K$ is a cable $C_{r_1, s_1}(T_{a, b}), |a|>|b|>1$.
    \end{enumerate}
\end{thm} 
\begin{proof}
We first show that Proposition \ref{prop:sxpiece} is realized. Suppose $K'$ is a knot such that there exists an orientation-preserving homeomorphism $S^3_{K}(p/q) \cong S^3_{K'}(p/q)$. If $S^3_{K}(p/q)$ is Seifert fibred, Proposition \ref{prop:sxpiece} is immediately realized. 

If $S^3_{K}(p/q)$ contains a JSJ torus, then the surgered piece $X'$ of $S^3_{K'}(p/q)$ is not a filling of a knot exterior. It follows that the JSJ decomposition of $S^3_{K'}$ does not contain hyperbolic pieces if $|q|>2$. The surgered piece $X'$ is thus Seifert fibred and it is a filling of a JSJ piece $Y'$ of $S^3_{K'}$ that is either a composing space or a cable space (Proposition \ref{prop:geomsxpiece}). 

If $Y'$ is a composing space, then by Lemma \ref{lemma:compositeknotext} and Section \ref{sec:scenario}, $X'$ must be the image by the homeomorphism of the surgered piece of $S^3_{K}(p/q)$ if $|q|>2$.

If $Y'$ is a cable space and if $X'$ is not the image of the surgered piece of $S^3_{K}(p/q)$, then $X'$ is the image of the exterior of a torus knot $T_{a,b}$ in $S^3_K$ (Proposition \ref{prop:geomswitch}). Using the notation introduced in Section \ref{sec:scenario}, let $\mathcal{T}= \del S^3_{T_{a,b}}$ and let $P$ be the pattern such that $\mathcal{T}$ decomposes $K$ into $P$ and $T_{a,b}$. 

Suppose $K$ is not an $n$-times iterated cable of $C_{r_1, s_1}(T_{a, b}), n \geq 1$. Then the pattern space $V_P$ contains a composing space that shares all its boundary components with other JSJ pieces of $S^3_K$. By the proof of Proposition \ref{prop:notSeifert fibre space}, $V_P(\mathcal{T}; q(w')^2/y)$ is not Seifert fibred. Applying the Cyclic surgery theorem (Theorem \ref{thm:cyclic}) as described in Section \ref{sec:sxpiece}, we obtain that $|q|=1$, contradicting (i). 

If $K$ is an $n$-times iterated cable of $C_{r_1, s_1}(T_{a, b}), n \geq 1$, we apply the same method as in Section \ref{sec:surgeredcable} to compare the distances between the regular fibre slope and the longitudinal slope on $\mathcal{T}$ and $\mathcal{T}'$. This yields the inequality $|q| \leq |r_1|+|a|$, which implies that if (ii) holds, $X'$ must be the image of the surgered piece of $S^3_{K}(p/q)$.

The theorem now follows, as (i), (ii) and (iii) are greater than or equal to the bounds from Section \ref{sec:proof}.
\end{proof}

Note that Theorem \ref{thm:SFSbounds}(i) is equivalent to Theorem \ref{thm:cableSFS} when applied to cable knots. If $K$ is not a cable knot in Theorem \ref{thm:cableSFS}, then it is a composite knot and the result follows from Theorem \ref{thm:composite}, which we prove in the next section.

\begin{thm:cableSFS}
    If $K$ is a knot with an exterior consisting solely of Seifert fibred JSJ pieces, with one of them being a composing space, then any slope that is neither integral nor half-integral is a characterizing slope for $K$.
\end{thm:cableSFS}

\begin{rem}
We relied on the constructive nature of Seifert fibred spaces to compute the above bounds. If $S^3_K$ contains hyperbolic JSJ pieces, the task becomes more difficult. Indeed, for generic cases, we need to determine values that realize Theorems \ref{thm:finitehyperb} and \ref{thm:CKhyperb}. Recently, Wakelin established a lower bound on $|q|$ for a slope $p/q$ to be characterizing for certain hyperbolic patterns (\cite{wakelin}). In forthcoming work with Wakelin, we combine her findings and our study of Seifert fibred JSJ pieces to obtain further results.
\end{rem}

  \section{Characterizing slopes for composite knots}\label{sec:theor2}

We now turn to the proof of Theorem \ref{thm:composite}. 

\begin{thm:composite}
    If $K$ is a composite knot, then every non-integral slope is characterizing for $K$.
\end{thm:composite}

\subsection{The surgered submanifold}
Thus far, we have made the assumption $|q|>2$, allowing us to define the surgered piece of a surgery along a knot. When $|q|=2$, the surgered piece may not be defined as in Definition \ref{def:sxpiece} if the resulting manifold is obtained from filling a hyperbolic JSJ piece. Indeed, the surgery operation can create essential tori, or it might yield a Seifert fibre space which admits a Seifert fibred structure that extends to other JSJ pieces. 

\begin{defn}\label{defn:surgeredmanifold}
Let $Y_0 \cup Y_1 \cup Y_2 \cup \ldots \cup Y_n$ be the JSJ decomposition of the exterior of a knot $K$.

If $|q|>1$, then up to re-indexing the $Y_i$, the JSJ decomposition of $S^3_K(p/q)$ is of the form
\[(X_0 \cup X_1 \cup \ldots \cup X_m) \cup (Y_i \cup  Y_{i+1} \cup \ldots \cup Y_n),\]
for some $1 \leq i \leq n$ and $m \geq 0$, and where none of the $X_j$'s are JSJ pieces of $S^3_K$. The manifold $X_0 \cup X_1 \cup \ldots \cup X_m$ is the \textit{surgered submanifold} of $S^3_K(p/q)$.

If the surgered submanifold is a JSJ piece of $S^3_K(p/q)$, i.e., $m=0$, then we may also call it the \textit{surgered piece} of $S^3_K(p/q)$.
\end{defn}

\begin{rem}\label{rem:hypfill}
This definition is compatible with Definition \ref{def:sxpiece}. In fact, the surgered submanifold of a surgery $S^3_{K'}(p/q)$ may not be a surgered piece only in the case where the outermost piece of $S^3_{K'}$ is hyperbolic (proof of Proposition \ref{prop:JSJsurgery} and Proposition \ref{prop:cablepqs2}).
\end{rem}

We obtain an analogue of Proposition \ref{prop:sxpiece} for composite knots when $|q|>1$.

\begin{prop}\label{prop:compositesxpieces}
Let $K$ be a composite knot and suppose there exists an orientation-preserving homeomorphism $S^3_{K}(p/q) \cong S^3_{K'}(p/q)$ for some knot $K'$.
If $|q|>1$, then the homeomorphism $S^3_{K}(p/q) \cong S^3_{K'}(p/q)$ carries the surgered piece of $S^3_K(p/q)$ to a JSJ piece of the surgered submanifold of $S^3_{K'}$.
\end{prop}
\begin{proof}
Let $K = K_1 \# K_2 \# \ldots \# K_n$ where the $K_i$'s are prime for each $i=1, \ldots, n$. Let $Y$ be the outermost composing space of $S^3_K$. It is homeomorphic to the exterior of the link in $S^3$ with unknotted components $L_0, L_1, \ldots, L_n$ such that each pair $(L_0, L_i)$ for $i = 1, \ldots, n$, is a Hopf link and the link formed by $L_1, \ldots, L_n$ is the unlink with $n$ components. Let $\mathcal{L}_i = \del \nu L_i$ be the boundary components of $Y$ (Figure \ref{fig:composing}). By Remark \ref{rem:hypfill}, the surgered piece $X = Y(\mathcal{L}_0; p/q)$ of $S^3_K(p/q)$ is well-defined.

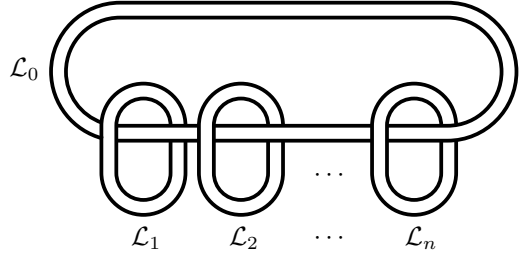
\begin{figure}[ht]
    \centering
    \begin{tikzpicture}
    \node at (0,0) {\includegraphics[scale=0.5]{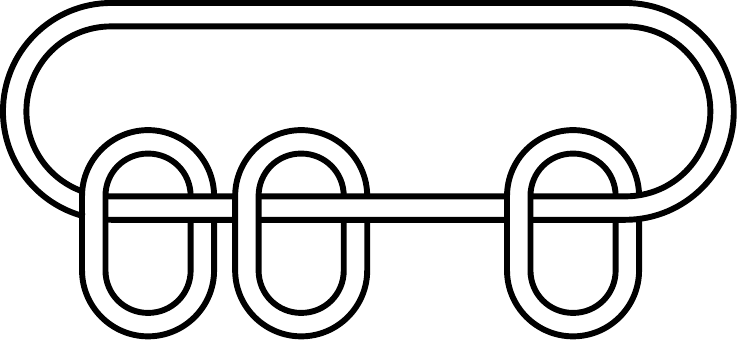}};

    \node at (-98pt,15pt) {$\mathcal{L}_0$};

    \node at (-52pt,-49pt) {$\mathcal{L}_1$};    
    \node at (-15pt,-49pt) {$\mathcal{L}_2$};    
    \node at (18pt,-49pt) {$\cdots$};    
    \node at (52pt,-49pt) {$\mathcal{L}_n$};   

    \node at (18pt,-25pt) {$\cdots$};    
    
    \end{tikzpicture}
    
    \caption{Composing space seen as a link complement in $S^3$}
    \label{fig:composing}
\end{figure}

Suppose the homeomorphism $S^3_{K}(p/q) \cong S^3_{K'}(p/q)$ does not carry $X$ into the surgered submanifold $X'$ of $S^3_{K'}$. Then there is a component $K_i$ of $K, i\in \{1, \ldots, n\}$, whose exterior contains a submanifold homeomorphic to $X'$. Let $\mathcal{T} = \del S^3_{K_i} \subset S^3_K$. The JSJ torus $\mathcal{T}$ decomposes $K$ into a composing pattern $P$ and $K_i$. The homeomorphism sends $\mathcal{T}$ to a JSJ torus $\mathcal{T}'$ of $S^3_{K'}(p/q)$ that separates $S^3_{K'}(p/q)$ into a manifold homeomorphic to $S^3_{K_i}$ and the exterior of a knot $J'$. As a result, $V_P(\mathcal{P};p/q)$ is homeomorphic to the exterior of $J'$, which contradicts Lemma \ref{lemma:compositeknotext}.
\end{proof}

\subsection{Fillings of a hyperbolic piece}\label{sec:fillingshyppiece}

In order to prove Theorem \ref{thm:composite}, we must demonstrate that the surgered submanifold of $S^3_{K'}(p/q)$ does not result from filling a hyperbolic JSJ piece of $S^3_{K'}$. Therefore, we now focus on the topology of the surgered submanifold of $S^3_{K'}(p/q)$, under the assumption that the outermost piece of $S^3_{K'}$ is hyperbolic. In this subsection, we study the surgery $S^3_{K'}(p/q)$ by itself, without taking into account any constraints arising from an orientation-preserving homeomorphism $S^3_{K}(p/q) \cong S^3_{K'}(p/q)$.

Recall from Theorem \ref{thm:budney} that if $Y'$ is a hyperbolic JSJ piece in the exterior of a knot, then $Y'$ is either the exterior of a hyperbolic knot $L_0'$ in $S^3$ or the exterior of a hyperbolic link in $S^3$ with components $L_0', L_1', \ldots, L_n'$ such that the link formed by $L_1', \ldots, L_n'$ is the unlink with $n$ components. From now on, we will denote the boundary components of such a hyperbolic piece $Y'$ by $\mathcal{L}_i' = \del \nu L_i'$, $i=0, \ldots, n$.

\begin{prop}\label{prop:notfillingofcomposite}
Let $Y'$ be a hyperbolic JSJ piece of a knot exterior. If $Y'(\mathcal{L}_0';p'/q')$ is homeomorphic to either a composing space or to a $p/q$-filling of a composing space where $|q| > 1$, then $|q'| \leq 1$.
\end{prop}
\begin{proof}
Let $Y$ be a composing space with $n+1$ boundary components labelled as in Figure \ref{fig:composing}. Whitout loss of generality, suppose that $Y'(\mathcal{L}_0';p'/q')$ is homeomorphic to $Y(\mathcal{L}_0;p/q)$, a Seifert fibre space with one exceptional fibre of order $|q|>1$. Then $Y'$ also has $n+1$ boundary components. Up to permuting indices, we can assume that for each $i=1, \ldots, n$, the homeomorphism maps $\mathcal{L}_i'$ to $\mathcal{L}_i$.

There exists infinitely many slopes $\sigma_1$ on $\mathcal{L}_1'$ such that the cores of the surgery solid torus is an exceptional fibre in $Y'(\mathcal{L}_0', \mathcal{L}_1'; p'/q', \sigma_1)$. There also exists infinitely many slopes $\sigma_i$ on each $\mathcal{L}_i'$, $i=2, \ldots, n$, such that the core of the surgery solid tori are regular fibres in $Y'(\mathcal{L}_0', \mathcal{L}_i'; p'/q', \sigma_i)$. Since hyperbolic manifolds possess only finitely many exceptional surgery slopes on each of their torus boundary components, we can choose $\sigma_1, \ldots, \sigma_n$ such that $\widetilde{Y'} = Y'(\mathcal{L}_1', \ldots, \mathcal{L}_n' ; \sigma_1, \ldots, \sigma_n)$ is hyperbolic. 

Now,  $\widetilde{Y'}(\mathcal{L}_0'; p'/q')$ has base orbifold $S^2$ with two exceptional fibres, which means that it has cyclic fundamental group. On the other hand, $Y'(\mathcal{L}_0'; 1/0)$ is homeomorphic to the exterior of the unlink with $n$ components. Therefore, $\widetilde{Y'}(\mathcal{L}_0'; 1/0)$ is a connected sum of manifolds with cyclic fundamental groups. By Boyer and Zhang (\cite[Corollary 1.4]{boyerzhang}), or the Cyclic surgery theorem (Theorem \ref{thm:cyclic}) if the connected sum is trivial, we must have $|q'| = \Delta(p'/q',1/0) \leq 1$ since $\widetilde{Y'}$ is not Seifert fibred.

Suppose now that $Y'(\mathcal{L}_0';p'/q')$ is homeomorphic to a composing space with $n$ boundary components. The argument is similar to that above. There are infinitely many slopes $\sigma_i$ on each component $\mathcal{L}_i'$ of $\del Y'(\mathcal{L}_0';p'/q')$ such that the cores of the surgery solid tori corresponding to $\sigma_1, \sigma_2$ are exceptional fibres and the cores of the surgery solid tori corresponding to $\sigma_3, \ldots, \sigma_n$ are regular fibres. We can choose the $\sigma_i$'s so that $\widetilde{Y'} = Y'(\mathcal{L}_1', \ldots, \mathcal{L}_n' ; \sigma_1, \ldots, \sigma_n)$ is hyperbolic. Then $\widetilde{Y'}(\mathcal{L}_0'; p'/q')$ and $\widetilde{Y'}(\mathcal{L}_0'; 1/0)$ are fillings of a hyperbolic manifold that are respectively a manifold with cyclic fundamental group and a connected sum of manifolds with cyclic fundamental groups. We conclude as before with \cite[Corollary 1.4]{boyerzhang} or Theorem \ref{thm:cyclic}.
\end{proof}

\begin{lemma}\label{lemma:hyperbJSJuntouched}
Let $K'$ be a knot such that the outermost piece $Y'$ of $S^3_{K'}$ is hyperbolic. If $|q|>1$, then the JSJ tori of $S^3_{K'}$ are incompressible in $S^3_{K'}(p/q)$.
\end{lemma}
\begin{proof} We may assume that $K'$ is a satellite knot as otherwise, the statement is vacuously true. The surgered submanifold of $S^3_{K'}(p/q)$ contains $Y'(\mathcal{L}_0';p/q)$. Suppose $Y'(\mathcal{L}_0';p/q)$ has compressible boundary. As $S^3_{K'}(p/q)$ is irreducible, $Y'(\mathcal{L}_0';p/q)$ is also irreducible, so it must be a solid torus. This implies that $Y'$ has two boundary components and, therefore, as $Y'(\mathcal{L}_0'; 1/0)$ is also a solid torus, a result of Wu (\cite[Theorem 1]{wu_incompressibility_1992}) implies that $|q| = \Delta(p/q,1/0)\leq 1$, which contradicts our assumption $|q|>1$. Hence, $Y'(\mathcal{L}_0';p/q)$ has incompressible boundary and as a consequence, the JSJ tori of $S^3_{K'}$ are incompressible in $S^3_{K'}(p/q)$. 
\end{proof}

\begin{cor}\label{cor:hyperbSurgmanUnionofJSJ} Let $K'$ be a knot such that the outermost piece $Y'$ of $S^3_{K'}$ is hyperbolic. If $|q|>1$, then the surgered submanifold of $S^3_{K'}(p/q)$ is either $Y'(\mathcal{L}_0';p/q)$, or the union of $Y'(\mathcal{L}_0';p/q)$ and some Seifert fibred JSJ pieces of $S^3_{K'}$ sharing a boundary component with $Y'$ in $S^3_{K'}$.\qed
\end{cor}

\subsection{Non-integral toroidal surgeries}\label{sec:nonintegraltoroidalsx}
We now study the surgered submanifold of $S^3_{K'}(p/q)$ in the context of an orientation-preserving homeomorphism $S^3_{K}(p/q) \cong S^3_{K'}(p/q)$ where $K$ is a composite knot.

Recall that a 3-manifold is said to be \textit{toroidal} if it contains an essential torus, and \textit{atoroidal} otherwise. The following proposition narrows down our investigation to hyperbolic manifolds that admit a non-integral toroidal surgery.

\begin{prop}\label{prop:hyperbSurgmanToroidal} Let $K$ be a composite knot. Suppose there exists an orientation-preserving homeomorphism $S^3_{K}(p/q) \cong S^3_{K'}(p/q)$ where $K'$ is such that the outermost piece $Y'$ of $S^3_{K'}$ is hyperbolic.  If $|q|>1$, then $Y'(\mathcal{L}_0';p/q)$ is toroidal.
\end{prop}
\begin{proof} By Proposition \ref{prop:compositesxpieces}, the homeomorphism $S^3_{K}(p/q) \cong S^3_{K'}(p/q)$ sends the surgered piece $X$ of $S^3_K(p/q)$ to a JSJ piece of the surgered submanifold $X'$ of $S^3_{K'}(p/q)$.

Suppose $Y'(\mathcal{L}_0';p/q)$ is atoroidal. Then by Corollary \ref{cor:hyperbSurgmanUnionofJSJ}, $X'$ has trivial JSJ decomposition, which means that it is homeomorphic to $X$, a filling of a composing space. Hence, $Y'(\mathcal{L}_0';p/q)$ is homeomorphic to a submanifold of a filling of a composing space. 

Since $Y'(\mathcal{L}_0';p/q)$ has incompressible torus boundary components, it must be homeomorphic to either a filling of a composing space or a composing space, since these are the only submanifolds of $X'$ that have such a boundary. However, this contradicts Proposition \ref{prop:notfillingofcomposite}.
\end{proof}


\begin{lemma}
Let $Y'$ be a hyperbolic JSJ piece of a knot exterior with at least three boundary components. If $|q|>1$, then $Y'(\mathcal{L}_0';p/q)$ is atoroidal.
\end{lemma}
\begin{proof}
The filling $Y'(\mathcal{L}_0'; 1/0)$ is the complement of the unlink and it has compressible boundary. If $Y'(\mathcal{L}_0'; p/q)$ contains an essential torus, then by a result of Wu (\cite[Theorem 4.1]{wu_sutured_1998}), we have $|q| = \Delta(p/q, 1/0) \leq 1$, contradicting the assumption $|q|>1$.
\end{proof}

Eudave-Mu\~noz contructed in \cite{eudavemunoz} a family of hyperbolic knots that admit half-integral toroidal surgeries. These surgeries produce a union of two Seifert fibre spaces. Gordon and Luecke proved that if a hyperbolic knot admits a non-integral toroidal surgery, then it belongs to Eudave-Mu\~noz's family and the surgery slope is half-integral (\cite{luecke_non-integral_2004}).

\begin{lemma}
Let $K$ be a composite knot. Suppose there exists an orientation-preserving homeomorphism $S^3_{K}(p/q) \cong S^3_{K'}(p/q)$ where $K'$ is such that the outermost piece of $S^3_{K'}$ is hyperbolic. If $|q|>1$, then $K'$ is not a hyperbolic knot.
\end{lemma}
\begin{proof} 
    By Eudave-Mu\~noz, Gordon and Luecke, any non-integral surgery along a hyperbolic knot contains at most one essential torus. However, the surgery $S^3_{K}(p/q)$ contains at least two essential tori, the boundary components of the surgered piece being such tori.
\end{proof}

We obtain the next corollary by combining Proposition \ref{prop:hyperbSurgmanToroidal} and the two preceding lemmas.

\begin{cor}\label{cor:hyponethreeboundary}
Let $K$ be a composite knot. Suppose there exists an orientation-preserving homeomorphism $S^3_{K}(p/q) \cong S^3_{K'}(p/q)$ where $K'$ is such that the outermost piece $Y'$ of $S^3_{K'}$ is hyperbolic. If $|q|>1$, then $Y'$ has exactly two boundary components. \qed
\end{cor}

\begin{figure}[h]
    \centering
    \begin{tikzpicture}
    \node at (0,0) {\includegraphics[scale=0.48]{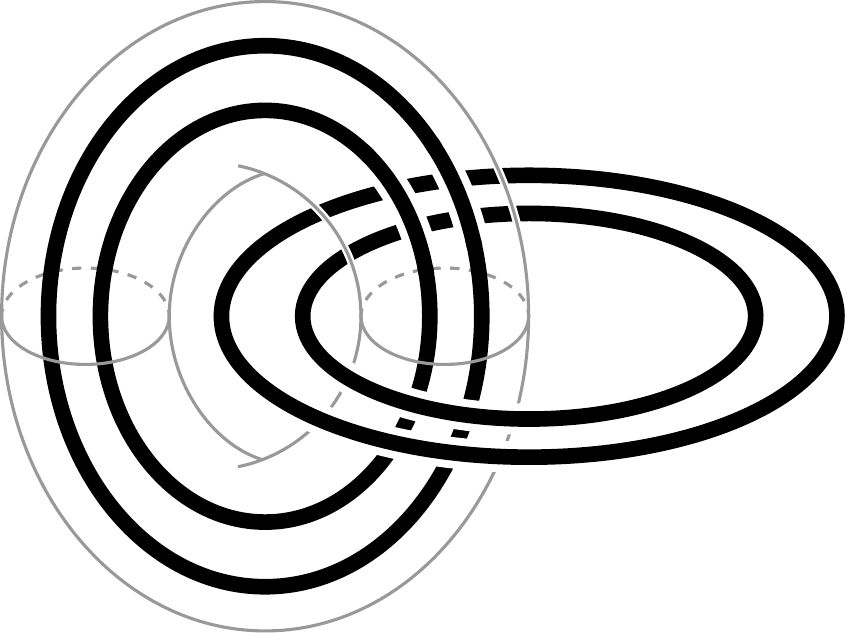}};

    \node at (-38pt,-37pt) {$(L_1, \alpha_1)$};    
    \node at (-38pt,-54pt) {$(L_2, \alpha_2)$};    

    \node at (55pt,0) {$(L_3, \alpha_3)$};   
    \node at (112pt,0) {$(L_4, *)$};   
    
    \end{tikzpicture}
    
    \caption{Surgery description of a non-integral toroidal surgery along a hyperbolic knot in $S^1 \times D^2$}
    \label{fig:toroidalsx}
\end{figure}
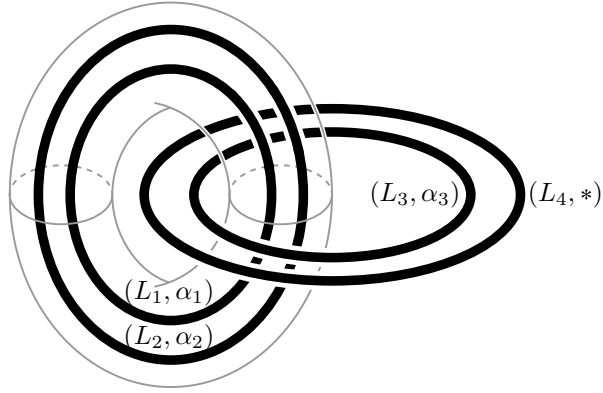

Gordon and Luecke also classified in \cite{luecke_non-integral_2004} all hyperbolic knots in solid tori which admit non-integral toroidal surgeries. They are derived from Eudave-Mu\~noz's construction mentioned above, and the resulting surgeries have link surgery descriptions as in Figure \ref{fig:toroidalsx}. The labels $L_i$ identify the link components and the $\alpha_i$'s are the corresponding surgery slopes, which correspond to the slopes $\alpha, \beta, \gamma$ in \cite{luecke_non-integral_2004}. The component $L_4$ is left unfilled. The essential torus $\widehat{\mathcal{T}}$ is pictured. If $\mu_i$ is the slope on $\del\nu L_i$ that bounds a disc in $\nu L_i$, then $\Delta(\alpha_i, \mu_i) \geq 2$ (\cite[Proof of Corollary A.2]{luecke_non-integral_2004}). Hence, the surgery is the union along $\widehat{\mathcal{T}}$ of two Seifert fibre spaces $M_1$ and $M_2$, with respective base orbifolds a disc with two cone points of orders $\Delta(\alpha_1, \mu_1)$ and $\Delta(\alpha_2, \mu_2)$, and an annulus with one cone point of order $\Delta(\alpha_3, \mu_3)$.

\begin{prop}[{\cite[Proof of claim in proof of Corollary A.2]{luecke_non-integral_2004}}]\label{prop:toroidalsxRegfibre}
    Let $\mathcal{E}$ be the exterior of a hyperbolic knot $K_0$ in $S^1 \times D^2$ such that $\mathcal{E}(\mathcal{K}_0; \sigma)$ is toroidal and $\Delta(\sigma, \mu) > 1$, where $\mu$ bounds a disc in $\nu K_0$.  Then $\mathcal{E}(\mathcal{K}_0; \sigma)$ is the union of Seifert fibre spaces $M_1$ and $M_2$. Suppose $\del M_2 = \del \mathcal{E}(\mathcal{K}_0; \sigma) = \del (S^1 \times D^2)$. The slope of a regular fibre of $M_2$ on $\del (S^1 \times D^2)$ does not coincide with the slope that bounds a disc in $\mathcal{E}(\mathcal{K}_0; \mu) \cong S^1 \times D^2$.
\end{prop}
\begin{proof}
We follow \cite{luecke_non-integral_2004}, in which the solid torus containing $K_0$ is denoted $L(\alpha, \beta, \gamma, *, 1/2)$ and the non-integral toroidal filling $\mathcal{E}(\mathcal{K}_0; \sigma)$ is denoted $L(\alpha, \beta, \gamma, *, 1/0)$. If $\Delta(\sigma, \mu)>1$ and $\mathcal{E}(\mathcal{K}_0; \sigma)$ is toroidal, then by the discussion above, $\mathcal{E}(\mathcal{K}_0; \sigma)$ is the union of Seifert fibre spaces $M_1$ and $M_2$, and its surgery description is given by Figure \ref{fig:toroidalsx}.

Let $h$ be the slope of a regular fibre of $M_2$ on $\mathcal{S} = \del (S^1 \times D^2)$. Suppose by contradiction that $h$ bounds a disc in $\mathcal{E}(\mathcal{K}_0; \mu) \cong S^1 \times D^2$. Then $\mathcal{E}(\mathcal{S}, \mathcal{K}_0; h, \mu) \cong S^2 \times S^1$. 

In the surgery description from Figure \ref{fig:toroidalsx}, filling along $\mathcal{S}$ corresponds to filling along $L_4$. One can see that filling $M_2$ along a regular fibre yields the connected sum of a lens space and a solid torus whose meridian has distance one with a regular fibre of $M_1$. Hence, $\mathcal{E}(\mathcal{S}, \mathcal{K}_0; h, \sigma)$ is either a lens space or a connected sum of lens spaces. Since $\Delta(\sigma, \mu)>1$, this implies that $\mathcal{E}(\mathcal{S}; h)$ is reducible (\cite[Theorem 1.2]{gordlue_reducible} and \cite[Corollary 1.4]{boyerzhang}). 

Write $\mathcal{E}(\mathcal{S}; h) = N_1 \# N_2$. Then $S^2 \times S^1 \cong \mathcal{E}(\mathcal{S}, \mathcal{K}_0; h, \mu) \cong N_1 \# N_2(\mathcal{K}_0; \mu)$. But $S^2 \times S^1$ does not contain a separating $S^2$, which means that $N_1 \cong S^2 \times S^1$ and $N_2(\mathcal{K}_0; \mu) \cong S^3$. It follows that $\mathcal{E}(\mathcal{S}, \mathcal{K}_0; h, \sigma) = (S^2 \times S^1) \# N_2(\mathcal{K}_0; \sigma)$. This contradicts the assertion that $\mathcal{E}(\mathcal{S}, \mathcal{K}_0; h, \sigma)$ is a lens space or a connected sum of lens spaces, as such spaces do not contain a non-separating essential sphere.
\end{proof}

\begin{prop}\label{prop:hyptwoboundary}
Let $K$ be a composite knot. Let $K'$ be a knot such that the outermost piece of $S^3_{K'}$ is hyperbolic. If $|q|>1$, then there is no orientation-preserving homeomorphism between $S^3_{K}(p/q)$ and $S^3_{K'}(p/q)$.
\end{prop}
\begin{proof}
    Suppose by contradiction that there exists an orientation-preserving homeomorphism \linebreak $S^3_{K}(p/q) \cong S^3_{K'}(p/q)$ where $K'$ is such that the outermost piece $Y'$ of $S^3_{K'}$ is hyperbolic, and where $|q|>1$. By Corollary \ref{cor:hyperbSurgmanUnionofJSJ}, the surgered submanifold of $S^3_{K'}(p/q)$ contains $Y'(\mathcal{L}_0';p/q)$. By Corollary \ref{cor:hyponethreeboundary} and Proposition \ref{prop:hyperbSurgmanToroidal}, $Y'$ is the exterior of a knot in a solid torus and $Y'(\mathcal{L}_0';p/q)$ is a non-integral toroidal filling. By Gordon and Luecke, $Y'(\mathcal{L}_0';p/q)$ is the union of two Seifert fibre spaces $M_1$ and $M_2$, with respective base orbifolds a disc with two cone points and an annulus with one cone point. 
    
    Let $X$ be the surgered piece of $S^3_{K}(p/q)$ and let $X'$ be its image in $S^3_{K'}(p/q)$ by the homeomorphism $S^3_{K}(p/q) \cong S^3_{K'}(p/q)$. Since $X'$ is a filing of a composing space, we have $X' \cap Y'(\mathcal{L}_0';p/q) = M_2$ (proof of Proposition \ref{prop:fillingY0}(2)).
    Let $\mathcal{T}' = \del Y'(\mathcal{L}_0';p/q) \subset \del M_2$. In $S^3_{K'}$, the torus $\mathcal{T}'$ decomposes $K'$ into $P'$ and $J'$. Let $\mathcal{E} = V_{P'}$ and $\mathcal{P}' = \del \nu P'$.
    
    The torus $\mathcal{T}'$ is the image by the homeomorphism of an incompressible torus $\mathcal{T}$ in $X \subset S^3_{K}$. Although this torus might not be a JSJ torus of $S^3_{K}$, it separates $S^3_K$ into a pattern space $V_P$ and a knot exterior $S^3_J$, where $P$ is a composing pattern (and $J$ is a composite knot if $\mathcal{T}$ is not a JSJ torus). Let $\mathcal{P} = \del \nu P \subset \del V_P$.

    We thus have homeomorphisms $V_P(\mathcal{P};p/q) \cong \mathcal{E}(\mathcal{P}';p/q)$ and $S^3_J \cong S^3_{J'}$. These imply that the meridian on $\mathcal{T}$ given by $J$ is sent by the homeomorphism $S^3_{K}(p/q) \cong S^3_{K'}(p/q)$ to the meridian on $\mathcal{T}'$ given by $J'$, by the knot complement theorem.

    By construction of satellite knots, the meridian on $\mathcal{T}'$ given by $J'$ coincides with the slope that bounds a disc in $\mathcal{E}(\mathcal{P}';1/0)$. By Proposition \ref{prop:toroidalsxRegfibre}, this slope does not coincide with the slope of a regular fibre of $M_2$ on $\mathcal{T}'$.
    On the other hand, a regular fibre on $\mathcal{T}$ in $V_P(\mathcal{P};p/q)$ comes from the Seifert fibred structure of a composing space, so it has meridional slope.

    Hence, regular fibres on $\mathcal{T}$ are not mapped to regular fibres on $\mathcal{T}'$. This contradicts the unicity of the Seifert fibred structure on a Seifert fibre space with base orbifold an annulus and one cone point.
\end{proof}

\begin{proof}[Proof of Theorem \ref{thm:composite}]
    Let $K$ be a composite knot and suppose there is an orientation-preserving homeomorphism $S^3_{K}(p/q) \cong S^3_{K'}(p/q)$ for some knot $K'$, where $|q|>1$.
    According to Proposition \ref{prop:compositesxpieces}, the surgered piece of $S^3_{K}(p/q)$ is carried into the surgered submanifold of $S^3_{K'}(p/q)$. Proposition \ref{prop:hyptwoboundary} implies that the surgered submanifold of $S^3_{K'}(p/q)$ is a JSJ piece of $S^3_{K'}(p/q)$, and it is the $p/(qt^2)$-filling of a JSJ piece $Y'$ of $S^3_{K'}$ for some $t \geq 1$. 
    
    According to Proposition \ref{prop:compositesxpieces}, this filling $Y'(p/(qt^2))$ is homeomorphic to the surgered piece of $S^3_{K}(p/q)$, a filling of a composing space with one exceptional fibre of order $|q|$. By Proposition \ref{prop:notfillingofcomposite}, $Y'$ is Seifert fibred. Since $Y'(p/(qt^2))$ has at least two boundary components, $Y'$ is a composing space by Theorem \ref{thm:budney}. The exceptional fibre of $Y'(p/(qt^2))$ has order $|qt^2|=|q|$, so $t=1$ and $K'$ is not a cable. Therefore, we conclude that $K = K'$.
\end{proof}

\bibliographystyle{amsalpha} 
\bibliography{cas-refs}

\end{document}